\documentclass[trsc]{informs3}
\OneAndAHalfSpacedXI 
\usepackage{eucal}

\usepackage{graphicx,booktabs,longtable}
\usepackage{bm}
\usepackage{microtype}
\usepackage{url}
\usepackage{caption}
\usepackage{subcaption}

\usepackage{natbib}
 \bibpunct[, ]{(}{)}{,}{a}{}{,}%
 \def\bibfont{\small}%

\usepackage{lscape}

\usepackage{pdflscape}
\ACKname{Acknowledgements}
\usepackage{afterpage}

\usepackage{tikz}
\usetikzlibrary{shapes,snakes,shapes.geometric, backgrounds, patterns, scopes,patterns,intersections,calc}

\usepackage{pgfplots}
\pgfplotsset{compat=1.11}
\usepgfplotslibrary{fillbetween}
\usetikzlibrary{intersections}
\usetikzlibrary{patterns}

\pgfdeclarelayer{bg}
\pgfsetlayers{bg,main}

\pgfdeclarelayer{bg}
\pgfsetlayers{bg,main}

\usepackage[mode=buildnew]{standalone}

\usetikzlibrary{arrows}

\usepackage[inline]{enumitem}

\usetikzlibrary{positioning}
\usepackage{listings}
\usepackage{color}
\usepackage{bold-extra}
\usepackage[ruled,linesnumbered]{algorithm2e}

\newcommand{\policy}{CFA}

\TheoremsNumberedThrough     
\EquationsNumberedThrough    
\usepackage{epstopdf}

\begin{document}

\RUNAUTHOR{Schrotenboer et al.}
\RUNTITLE{Routing for Local Delivery Platforms}
\TITLE{Fighting the E-commerce Giants: Efficient Routing and Effective Consolidation for Local Delivery Platforms}

\ARTICLEAUTHORS{%
	\AUTHOR{Albert H.~Schrotenboer}
	\AFF{Operations, Planning, Accounting, and Control Group, School of Engineering, Eindhoven University of Technology\\ \EMAIL{a.h.schrotenboer@tue.nl}}
	\AUTHOR{Michiel A.~J.~uit het Broek, Paul Buijs}
	\AFF{Department of Operations, Faculty of Economics and Business, University of Groningen,\\ \EMAIL{a.j.uit.het.broek@rug.nl, p.buijs@rug.nl}}
	\AUTHOR{Marlin W.~Ulmer}
	\AFF{Chair of Management Science, Faculty of Economics, Otto Von Guericke Universität Magdeburg, \\ \EMAIL{marlin.ulmer@ovgu.de}}
}

\ABSTRACT{Local delivery platforms are collaborative undertakings where local stores offer instant-delivery to local customers ordering their products online. Offering such delivery services both cost-efficiently and reliably is one of the main challenges for local delivery platforms, as they face a complex, dynamic, stochastic dynamic pickup-and-delivery problem. Orders need to be consolidated to increase the efficiency of the delivery operations and thereby enable a high service guarantee towards the customer and stores. But, waiting for consolidation opportunities may jeopardize delivery service reliability in the future, and thus requires anticipating future demand. This paper introduces a generic approach to balance the consolidation potential and delivery urgency of orders. Specifically, it presents a newly developed parameterized Cost-Function Approximation (CFA) approach that modifies a set-packing formulation with two parameters. This CFA approach not only anticipates future demand but also utilizes column generation to search the large decision space related to pickup-and-delivery problems fast. Inspired by a motivating application in the city of Groningen, the Netherlands, numerical experiments show that our CFA approach strongly increases perceived customer satisfaction while lowering the total travel time of the vehicles compared to various benchmark policies. Furthermore, our CFA also reduces the percentage of late deliveries, and their lateness, to a minimum. Finally, our approach may assist managers in practice to manage the non-trivial balance between consolidation opportunity and delivery urgency.}

\KEYWORDS{Local Delivery Platforms; Dynamic Pickup and Delivery; Cost-Function Approximation; Column Generation; City Logistics; Last-Mile Delivery}

\maketitle 

\section{Introduction}

``See now. Buy now. How fast can I get it?'' \citep{deloitte2019}. The demand for instant-delivery services is rising fast. Large, global e-commerce companies, such as Amazon, and their logistics service providers have anticipated this trend---or indeed, were a driving force behind it---and operate a distribution network that is increasingly capable of delivering an enormous range of products in very short time. In the meanwhile, local stores struggle to keep up. Many have seen the number of in-store shoppers drop consistently and have not been able to compensate these sales via a stronger online presence. The Coronavirus pandemic and ensuing lock down measures only exacerbated this trend.

\begin{figure}
    \centering
    \includegraphics[scale = 0.52]{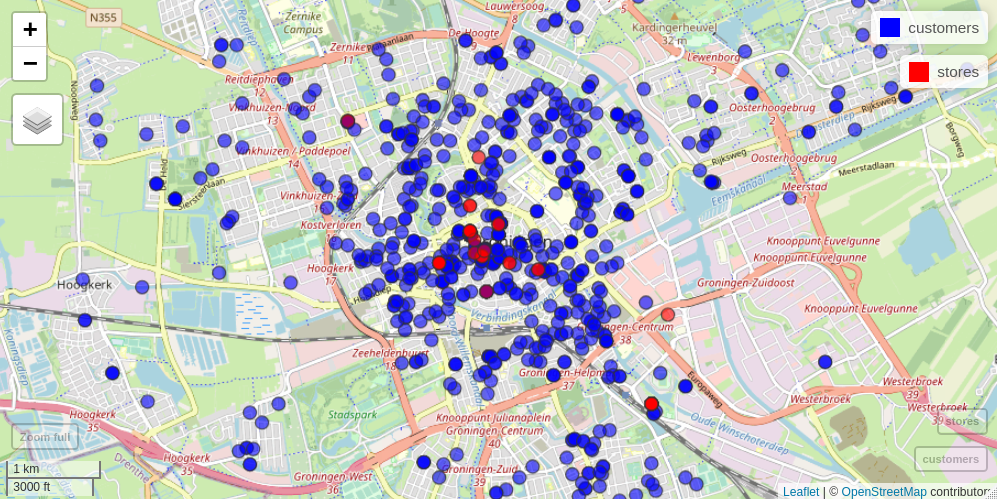}
        \caption{Order realizations of the days before Christmas 2020, with the stores in red and the customers in blue.}
    \label{fig:groningen}
\end{figure}

One way local stores can counter the enormous power of the e-commerce giants is by setting up a platform where local customers can order products from one or multiple stores and get those products delivered to their homes \citep{forbes}. Our study is inspired by an example from the Netherlands, where stores in the city of Groningen cooperate and, with the help of a local software developer and a local bicycle courier company, established and operate a local shopping platform called ``Warenhuis Groningen''. Customers can browse and order products at over 50 local stores. Requests for the pickup and delivery of the products purchased can be entered both by the customers and stores. Customers may require multiple products to be picked up from different stores for delivery at their home. Stores may require multiple products to be picked up at their store for delivery at different customers. This leads to a complex picture of pickups and deliveries spread across a city, as exemplified in Figure~\ref{fig:groningen}, which plots all orders of a single day shortly preceding Christmas 2020. On top of this already complex picture, the orders arrive dynamically throughout the day, and the platform must assign the orders expeditiously to couriers that operate a fleet of cargo-bikes to enable quick, reliable, and sustainable home delivery. While customers may tolerate occasional delays to support local stores, local delivery platforms face fierce competition from several e-commerce giants. Thus, these platforms must efficiently route and consolidate orders to stay alive in this competitive market environment.

In this paper, we address the dynamic, stochastic pickup-and-delivery problem faced by local delivery platforms, such as the one introduced above. One novelty of this complex optimization problem lies in the dynamic arrival of both one-to-$n$ and $n$-to-one orders (i.e., from a single store to multiple customers or from multiple stores to a single customer), as opposed to the typical e-commerce delivery models where many customers are served from one or multiple depots. Another important feature of the problem is that customers attribute a certain value to on-time delivery of their order and/or penalize varying degrees of late delivery. We model customer satisfaction as a general function---that can take any form---based on the time between ordering and delivery of the product and a and a predefined soft deadline. The problem's objective is to assign orders to cargo-bike couriers and find routes for those couriers in such a way that the expected customer satisfaction is maximized. We refer to this problem as Dynamic Pickup-and-Delivery Problem for Local Platforms (DPDP-LP).

Achieving high delivery service quality in the context of the DPDP-LP is challenging. The platform needs to balance \textit{the urgency of the orders currently known} with \textit{future consolidation opportunities} emerging from currently unknown future demand. For example, sending a courier instantly when an order arrives ensures fast delivery, but also neglects possibilities for consolidation in the future and consumes resources that may be required to deliver future demand on time. Further, assignment and routing decisions are very complex since a large-scale pickup-and-delivery problem with soft deadlines must be solved at every decision moment. The platform, therefore, requires advanced tools to optimize over a large decision space in reasonable computation times. To address this challenging problem, we develop a novel parameterized cost function approximation (CFA) approach. The idea of our CFA approach is to use the powerful mixed-integer programming tools from the Operations Research community to search the vast decision space in every state. To avoid myopic and inflexible solutions, our approach manipulates the mixed integer program (MIP) in its objective function and constraints. Such manipulations must be crafted carefully, and done with respect to domain knowledge \citep{powell2007approximate}.

Our CFA approach balances optimization and anticipation as follows: For searching the decision space and finding efficient order assignment and vehicle routing decisions in real-time, our approach uses a Dantzig-Wolfe reformulation of the pickup-and-delivery problem. The resulting set-packing formulation is solved using column generation techniques so that the large decision space can be searched efficiently. To ensure the decisions are also effective with respect to future demand, our approach modifies the set-packing formulation with two parameters. The first parameter, $\alpha$, considers ``how well an order fits into a vehicle path". The second parameter, $\beta$, accounts for ``how urgent a currently unscheduled order is". These parameters scale the reduced costs of the variables in the set-packing formulation, and by using column generation to generate new variables (i.e., vehicle paths), our CFA approach generates vehicle paths that automatically balances between urgency and consolidation opportunities. 

Our work makes several important contributions to the literature: 
\begin{itemize}
    \item We introduce a new and emerging problem in city logistics, the Dynamic Pickup-and-Delivery Problem for City Logistics (DPDP-LP). It models the fundamental problem of real-time fulfillment of $n$-to-1 and $1$-to-$n$ pickup-and-delivery requests with soft deadlines, with a variety of real-world applications in the context of city logistics. 
    \item We develop a novel CFA approach that combines advanced mixed-integer programming with anticipatory decision making. To the best of our knowledge, it is the first generic anticipatory approach that can handle exponentially large decision spaces for dynamic routing problems with a service-oriented objective. Our CFA approach relies on a parameterization of a set-packing model with only two parameters that consider how urgent orders are and what future consolidation opportunities can be expected. These parameters are easy to interpret and can be readily implemented in existing routing optimization tools used in practice. The best-performing parameter values can be easily found offline using global optimization tools. 
    \item We show the advantages of our CFA approach in a variety of systemically created instance settings. Our CFA approach both outperforms policies that focus on minimizing travel distances and policies that focus on instant delivery. Specifically, comparisons with a variety of benchmark policies reveal that our approach manages the time available before the soft order deadlines significantly better while also improving customer satisfaction. The best-performing parameter combinations in our CFA approach have a non-linear relation, suggesting that it is not trivial for managers to develop efficient policies without the support of our approach.
    \item Our CFA approach is applicable to a wide range of stochastic, dynamic pickup-and-delivery problems in that no structural adaptations are required to account for new business-specific conditions and constraints as is usually the case when relying on solutions based on meta-heuristic techniques.
    \item From a managerial point of view, we investigate the interplay between urgency and consolidation opportunities and derive easily interpretable decision support. Specifically, more emphasis should be placed on urgency when there are plenty vehicles available (i.e., many vehicles in relation to the number of outstanding orders) and when late delivery results in a fixed costs, regardless of the extent of the dela. In contrast, consolidation is relatively important when solely minimizing average order lateness, as typically the case in meal delivery services. Furthermore, we observe that our CFA approach is robust for different customer order sizes and the amount of stores appearing in the system. 
    \item Concluding, our CFA approach draws on advanced MIP-based optimization techniques and state-of-the-art cost-function approximation while not requiring enormous computational power involved with making new assignment and routing decisions at every newly arriving order. As such, it supports local delivery platforms in unleashing their full potential and compete against the e-commerce giants.  
\end{itemize}

The remainder of this paper is organized as follows. In Section~\ref{sec:litrev}, we position our work in the relevant literature. Section~\ref{sec:prob} presents the formal problem statement of the DPDP-LP, and Section~\ref{sec:mod} provides a Markov Decision Process Formulation. The Cost-Function Approximation (CFA) approach is described in Section~\ref{sec:sol}. Numerical results and managerial insights on the performance of our CFA approach and related policies are presented in Section~\ref{sec:result}. We conclude our paper in Section~\ref{sec:conc}

\section{Literature Review}
\label{sec:litrev}

The stochastic, dynamic pickup-and-delivery problem that arises at local delivery platforms, as studied in this paper, and the CFA that we propose for efficient routing and effective consolidation policies will be discussed from two perspectives. First, we consider the deterministic pickup-and-delivery literature and identify the state-of-the art approaches for efficient routing. Second, we outline the literature that relates to instant deliveries and focus on how dynamic routing is done in such applications. For the latter part, a concise overview of related work is provided in Table~\ref{tab:lit}. For reviews on dynamic pickup-and-delivery problems in general, we refer the reader to the review by \cite{berbeglia2010dynamic}.

\subsection{Deterministic Pickup-and-Delivery Routing}

Local delivery platforms need to assign pickup-and-delivery orders to a fleet of available vehicles in real-time to fulfill the associated customer demands. At a fixed decision epoch in time, this problem resembles a many-to-many pickup-and-delivery problem with delivery deadlines (see the classification by \cite{berbeglia2007static}). Our problem has two distinguishing elements: 1) the complex network structure arising from one-to-n and n-to-one pickup and delivery orders (from a single store to multiple customers and from multiple stores to a single customer, respectively), and 2) the objective that seeks to maximize customer satisfaction by minimizing a general \textit{penalty function} that maps the delivery time of an order to a cost in relation to a soft delivery deadline associated with the order.

From the viewpoint of each single product ordered, a clear one-to-one relation exists between the pickup at the store and the delivery at the customer, and our proposed solution approach routes each product independently. This bears resemblance with the static one-to-one pickup and delivery problem (see, \cite{parragh2008survey} and \cite{braekers2016vehicle}), for which efficient exact solution methods are readily available (see, \cite{ropke2006adaptive} and \cite{ropke2009branch}). These solution methods focus heavily on finding efficient consolidation opportunity across orders to minimize the distance of the typically vehicle routes, while adhering to the deadlines involved with those orders. Due to the static nature of the problem, there are often many opportunities for consolidation, leading to relatively long, efficient routes.

In practice, delivery deadlines are not always hard. Especially in settings with same-day delivery services, customers value quick delivery. Providers may offer a certain service promise (e.g., delivered within 2 or 4 hours after ordering) which they try to meet while maximizing customer satisfaction. The inclusion of soft delivery deadlines in deterministic pickup-and-delivery relates to the recent work by \cite{he2019branch}, in which a branch-and-price-and-cut algorithm is presented that deals with travel costs that are a convex function of the specific arrival time within pre-specified time windows. To the best of our knowledge, costs as a function of the time it takes to deliver an order in pickup-and-delivery problems have only be considered in the context of meal deliveries \cite{yildiz2019provably}, where the penalty function takes a linear form. 

We deem applying a complete branch-and-price-and-cut method to our problem context infeasible due to the (often) long computation times involved with such methods. However, the concept of column-generation \citep{lubbecke2005selected, desaulniers2006column} via heuristic procedures \cite[see, e.g., ][]{desaulniers2008tabu} is considered very efficient in achieving high quality solutions. For instance, optimality gaps less then~0.5\% are reported for heuristic column generation procedures in \cite{sonntagtactical}. Typically, solution methods only generate columns via a heuristic in the root node and, afterward, assign routes (i.e., columns) via a set-partitioning alike model. Our CFA builds upon this approach, but by parameterizing the a set-partitioning formulation, we directly generate columns that are both good in terms of their combinatorial solution (i.e., routing efficiency) and in terms of anticipation of future events and decisions.

\subsection{Stochastic and Dynamic Pickup-and-Delivery Routing}
\begin{table}[h]
\TABLE 
{Overview of related works \label{tab:lit}}
{\begin{tabular}{l|ccccc|ccc}
\toprule  

 & \multicolumn{5}{c}{Model features} & \multicolumn{3}{c}{Method features} \\ \cmidrule(lr){2-6} \cmidrule(lr){7-9}
& \rotatebox{80}{Large decision space}  & \rotatebox{80}{Service objective}   &\rotatebox{80}{Accept all customers} & \rotatebox{80}{Multiple Depots}   & \rotatebox{80}{Dispatching } &  \rotatebox{80}{Generic Routing} & \rotatebox{80}{Anticipation} & \rotatebox{80}{Offline Training} \\ \midrule
   \textbf{This paper}  & \textbf{\checkmark} & \textbf{\checkmark} & \textbf{\checkmark }& \textbf{\checkmark} & \textbf{\checkmark} & \textbf{\checkmark} & \textbf{\checkmark}  & \checkmark\\ \midrule
 \cite{ulmer2021restaurant}  &                   & \checkmark & \checkmark & \checkmark &\checkmark   & & \checkmark & \checkmark \\ 
 \cite{steever2019dynamic}    &     & \checkmark & \checkmark & \checkmark &\checkmark   & & \checkmark &  \\
 \cite{reyes2018meal}        & & \checkmark & \checkmark & \checkmark & \checkmark  & & &  \\ \midrule
 \cite{arslan2019crowdsourced}&      &  & \checkmark & \checkmark & \checkmark  & & &\\
\cite{arslan2019splitting} &  &        \checkmark  & \checkmark & \checkmark &\checkmark  & $\sim$   & &  \\ \midrule
\cite{voccia2017same}  & \checkmark       & \checkmark &  & && $\sim$  & \checkmark \\  
\cite{azi2012dynamic}   &  \checkmark        & \checkmark &  &  & &$\sim$  & \checkmark \\ \midrule
\cite{sheridan2013dynamic}& \checkmark  &    &  \checkmark & \checkmark &  &  & & \\
\cite{ghiani2009anticipatory}&\checkmark  & \checkmark & \checkmark  & \checkmark & &$\sim$  & $\sim$  \\
\cite{mitrovic2004double}  &\checkmark &   & \checkmark   & \checkmark & &$\sim$  & $\sim$  & \\ 
\cite{mitrovic2004waiting}& \checkmark  &    & \checkmark   & \checkmark & &$\sim$  & $\sim$  &  \\ 

 \bottomrule
\end{tabular}}
{}
\end{table}
Table~\ref{tab:lit} presents an overview of the problem characteristics and decisions incorporated in our work, compared to the literature on dynamic and stochastic pickup-and-delivery routing. In the columns, we denote major features of the models and methods relevant to our context. In what follows, we classify related work according to these features to clearly position our model and method.

A check-mark $(\checkmark)$ in the column \textit{Large decision space} indicate that the model involves complex vehicle routing decisions rather than assignment-based decisions. A check-mark $(\checkmark)$ in the column \textit{Service objective} implies that a model is aimed at maximizing an objective related to customer service as opposed to the more typical minimization of travel distance, costs, or time. The check-marks in column \textit{Accept all customers} indicate that all customers that appear in the system need to be served. Alternatively, when not checked, the model allows for rejecting some customers , usually at a given penalty cost. Models that maximize service without the need to accept all customers refers to settings where the goal is to maximize the number of customers that are served. The column \textit{Multiple Depots} is marked if the pickup of products is done at multiple locations. The final column in the model category indicates whether vehicles are \textit{dispatching} and no route diversion are allowed after assigning requests.

Regarding the method features, a check-mark in the column \textit{Generic Routing} implies that the method employs generic techniques that can handle an exponentially large decision space faced at each decision epoch. A `$\sim$'-symbol implies problem specific meta-heuristics are used, and no symbol indicates that the routing problem is either solved with simple heuristics or its decision space is so constrained that it can be enumerated or it is solved with relatively straightforward heuristic procedures. The column \textit{Anticipation} indicates whether anticipation, for example using concepts of approximate dynamic programming (ADP), are used in a method, where a `$\sim$'-symbol means some form of anticipation but not using ADP ideas. The final column \textit{Offline Training} indicates that the method can be trained offline. The papers that provide some form of anticipation but without offline sampling typically require computationally expensive online scenario sampling.

A first research stream of relevance is meal-delivery routing, which considers the stochastic, dynamic pickup-and-delivery routing to transport meals between restaurants and customers. Both \cite{reyes2018meal} and \cite{ulmer2021restaurant} adopt a service-oriented objective value, which is modeled in a linear fashion, and only by minimizing the exceedance of deadline. The DPDP-LP addressed in our paper generalizes these models by considering a general penalty function that maps the time needed for delivery to a cost, considering a soft deadline stemming from, for example, a service promise. This penalty function can take any shape, such as a function where cost increases exponentially after a soft deadline is exceeded or where exceeding a deadline results in a fixed cost. More importantly, the DPDP-LP departs from meal-delivery routing by considering general parcels or products instead of meals. This usually allows for way more opportunities for consolidating orders as in the delivery of general parcels concern for freshness and temperature restrictions is not needed. The large number of consolidation opportunities result in an exponentially large decision space for the vehicle routes that can be constructed at each decision epoch. In contrast, the problem settings considered in \cite{reyes2018meal} and \cite{ulmer2021restaurant} involve fairly restricted decision spaces, which allows for using fairly simple routing heuristics. Our CFA approach is able to achieve the same anticipation in assigning orders to vehicles as considered by the Anticipatory Customer Assignment method of \cite{ulmer2021restaurant}, yet approximates individual customer order urgency to a cost function generalizing upon the inclusion of a time buffer before the soft deadline as in \cite{ulmer2021restaurant}. 

A second research stream of literature that relates to the DPDP-LP focuses on crowd-sourced delivery systems \cite{arslan2019crowdsourced, arslan2019splitting}. In such problems, drivers are not managed by a central decision-maker but appear randomly in the system to perform pickup-and-delivery tasks. In \cite{arslan2019crowdsourced}, the products ordered from a single customer cannot be split among multiple vehicles, while in \cite{arslan2019crowdsourced} this is allowed. In the DPDP-LP, we take the latter perspective. In contrast to these papers, we consider a known fleet of vehicles---and their couriers---that are managed by a central entity. More importantly, the DPDP-LP departs from crowd-sourced delivery  systems by having a larger decision space, which limits the applicability of the methods proposed for crowd-sourced to the DPDP-LP. This has two main reasons. First, from a problems perspective, the DPDP-LP focus solely on minimizing the incurred penalty costs related to customer service and do not incorporate driver-specific costs in crowd-sourced delivery systems. Second, in crowd-sourced delivery systems there is less attention required for anticipating future demand since the problem is simultaneously constrained by a limited number of drivers that become available when they want and often for short duration. In the DPDP-LP, drivers are managed by a central entity, which constructs routes to fulfill pickup and deliveries throughout the planning horizon.

A third research stream of literature that relates to our work is that of same-day delivery services, in which a single location (e.g., a retail store) receives orders dynamically \citep{azi2012dynamic, voccia2017same}. The important decision in this context is to determine when a vehicle should be dispatched from the depot. Similar to the DPDP-LP and meal-delivery routing problems, postponing the dispatching decision increases the potential to consolidate orders at the expense of less resources and flexibility in the future. \cite{voccia2017same} aim to maximize the numbers of customers that can be accepted while allowing disregarding incoming orders at some fixed cost. In \citep{azi2012dynamic}, an Adaptive Large Neighborhood Search is presented in combination with sampling scenarios of future demand \citep[see, e.g.,][]{bent2004scenario}. The authors consider minimizing a weighted sum of customer profits (independent of arrival time) minus travel distances. In contrast, in the the DPDP-LP we consider multiple stores and a customer service maximization. 

A fourth research stream of literature addresses general dynamic pickup-and-delivery problems, typically with hard time windows. \cite{mitrovic2004double, mitrovic2004waiting} acknowledge the trade-off between routing optimization and the urgency of fulfilling a pickup-and-delivery request, but solely with the objective to minimize routing distance. \cite{ghiani2009anticipatory} consider service objectives, but their work focuses on scenario sampling in combination with a routing heuristic. All of these works also allow for en-route diversions assuming last-minute arrival time communication to the customers, see for instance \cite{sheridan2013dynamic}, which is not allowed in the DPDP-LP, where customers require arrival time information once a routing decision is made. The power of our CFA approach lies in its applicability in practice by only relying on offline learning of the parameters in the CFA. Thus, our method eradicates the need for computationally expensive need of sampling future scenarios in a real-time fashion.

\section{Model}
\label{sec:prob}


This section first gives a problem narrative of the DPDP-LP, introducing the global notation of the system. To illustrate the decision space of the DPDP-LP, we continue with an example that also shows how the different aspects of the DPDP-LP come into play. Then, we provide a formal model of the DPDP-LP in the form of a stochastic sequential decision process.

\subsection{Problem Narrative}

The DPDP-LP considers the delivery of products that are ordered on a local delivery platform. These delivery processes take place on the time horizon $\cal T$. Orders are placed within this time horizon but can be delivered after this horizon. Connected to the platform are a set of stores $\cal D = \{1, \ldots, D\}$ that are known at the start of the time horizon. Each store $d \in \cal D$ has a known geographical location at $(\textsc x^d, \textsc y^d)$. 

The platform faces two types of dynamically arriving stochastic pickup-and-delivery orders. The first type is associated with a single store sending products via the platform to one or multiple customers, referred to as 1-to-$n$ orders. The second type is associated with a single customer ordering products from one or multiple stores via the platform, referred to as $n$-to-1 orders. A customer location will typically be the customer's home, but could also be a parcel pickup-point, a self-service parcel locker, or any location where the customer wants the parcel to be delivered. Pickup-and-delivery orders are only observed upon ordering, including information on how many parcels are to be picked up from which store(s) and/or to be delivered to which customer(s). We refer to each individual product of an pickup-and-delivery order as a \textit{request}.

A set of vehicles $\cal V = \{1, \ldots, V\}$, cargo-bikes in our example, is available to fulfill the pickup-and-delivery requests. All vehicles start at coordinate $(x^0, y^0)$, are available throughout the entire time horizon, and are assumed to have unlimited capacity. Because the transportation requests are relatively urgent and there are often vehicles of different sizes available, the vehicle capacity is unlikely to be a constraining factor in many practical scenarios \citep{allen2018understanding, chiara2021understanding}. The travel time between stores and customers is given by the deterministic travel time function $d(\cdot, \cdot)$.

Over the course of the time horizon, pickup-and-delivery requests are assigned to vehicles. At each point in time, idling vehicles can be assigned to perform yet unassigned requests. The assignment of requests involves a partitioning of the unassigned request over the idle vehicles, and to determine a path for each vehicle connecting the store and customer locations of the assigned requests. We do not allow en-route diversions, that is, vehicle routes cannot be altered nor extended once started. Therefore, it is also not needed to consider busy vehicles when assigning requests. Not allowing for en-route diversions enables the platform to communicate estimated arrival times to the stores and customers at the moment a request is assigned to a vehicle, which is inspired by current practice at our motivating application. By extension, we assume that both the pickup and the delivery of a request are assigned simultaneously. Note, we do allow multiple products of a single order to be split across multiple vehicles, and we do allow multiple requests of different orders to be consolidated and transported by a single vehicle. 

The platform has a service promise indicated by a soft delivery deadline that is fixed and equal to~$\bar{d}$ time units after the request is revealed to the system. For instance, if a request is placed at time $t \in \cal T$, then the delivery deadline is $t + \bar{d}$.

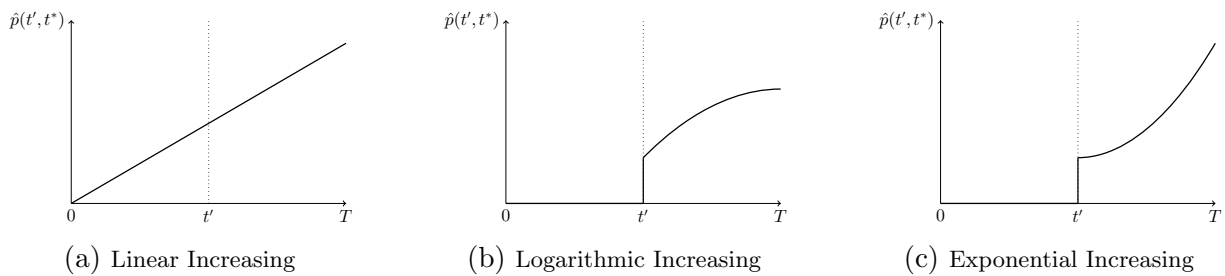
\begin{figure}
\centering
     \begin{subfigure}[b]{0.3\textwidth}
         \centering
         \resizebox{.95\linewidth}{!}{\begin{tikzpicture}
\draw[->] (0,0) -- (6,0) node[anchor=north] {$T$};
\draw	(0,0) node[anchor=north] {$0$}
		(3,0) node[anchor=north] {$t'$};

\draw[->] (0,0) -- (0,4) node[anchor=east] {$\hat p(t',t^*)$};
\draw[dotted] (3,0) -- (3,4);
\draw[thick] (0,0) -- (6,3.5);

\end{tikzpicture}}
         \caption{\footnotesize Linear Increasing}
         \label{fig:penLin}
     \end{subfigure}
     \hfill
     \begin{subfigure}[b]{0.3\textwidth}
         \centering
         \resizebox{.95\linewidth}{!}{\begin{tikzpicture}
\draw[->] (0,0) -- (6,0) node[anchor=north] {$T$};
\draw	(0,0) node[anchor=north] {$0$}
		(3,0) node[anchor=north] {$t'$};

\draw[->] (0,0) -- (0,4) node[anchor=east] {$\hat p(t', t^*)$};
\draw[dotted] (3,0) -- (3,4);

\draw[thick] (0,0) -- (3,0);
\draw[thick] (3,0) -- (3,1);
\draw[thick] (6,2.5) parabola (3,1);

\end{tikzpicture}}
         \caption{\footnotesize Logarithmic Increasing}
         \label{fig:penLog}
     \end{subfigure}
     \hfill
     \begin{subfigure}[b]{0.3\textwidth}
         \centering
         \resizebox{.95\linewidth}{!}{\begin{tikzpicture}
\draw[->] (0,0) -- (6,0) node[anchor=north] {$T$};
\draw	(0,0) node[anchor=north] {$0$}
		(3,0) node[anchor=north] {$t'$};

\draw[->] (0,0) -- (0,4) node[anchor=east] {$\hat p(t', t^*)$};
\draw[dotted] (3,0) -- (3,4);

\draw[thick] (0,0) -- (3,0);
\draw[thick] (3,0) -- (3,1);
\draw[thick] (3,1) parabola (6,3.5);

\end{tikzpicture}}
          \caption{\footnotesize Exponential Increasing}
         \label{fig:penExp}
     \end{subfigure}
    \caption{Three examples of penalty function $\hat p$}
    \label{fig:penalty_examples}
\end{figure}

The perceived customer satisfaction is determined by the time of delivery at the customer location. Let $\hat p(t', t^*)~:~[0, T]^2~\rightarrow~\mathbb{R}_{\geq 0}$ be a function that maps the delivery at time $t^* \in \cal T$ with associated delivery deadline $t'$ to a penalty $\hat{p}(t',t^*) \in \mathbb{R}_{\geq 0}$. Although this penalty function may take any form, we consider it to be strictly positive for $ t^* > t$. We interpret this penalty as a proxy for the (inverse) perceived customer satisfaction: the lower the penalty the more satisfied the customer is. Several examples that model different types of customer (dis)satisfaction are presented in Figure~\ref{fig:penalty_examples}. Figure~\ref{fig:penLin} shows a linear increasing penalty function, resembling the situation that requests should be delivered as-soon-as-possible; Figure~\ref{fig:penLog} a logarithmic increasing penalty function resembling that being late is a major problem and additional costs diminish for increasing deadline exceedances; and Figure~\ref{fig:penExp} shows an exponentially increasing penalty function resembling where customers are getting increasingly dissatisfied for for larger deadline exceedances. 


The goal of the DPDP-LP is to determine a policy that minimizes the expected incurred penalty costs, while all the stochastic and dynamic customer orders are fulfilled in real-time.


\begin{figure}[h!]
   \centering
   \raisebox{0.5\height}{\usetikzlibrary{shapes, snakes}
\usetikzlibrary{shapes.geometric}

\begin{tikzpicture}

\node[] at (0.4, 0.4) {Stores};
\node[rectangle, draw, fill, red, very thick] (r) at (0, 0) {};
\node[rectangle, draw, fill, blue, very thick] (b) at (0.4, 0) {};
\node[rectangle, draw, fill, brown, very thick] (br) at (0.8, 0) {};

\node[circle, draw, red, fill=red!40, very thick] (br) at (2.8,0) {};

\node[circle, draw, blue, fill=blue!40, very thick] (br) at (3.2, 0) {};

\node[circle, draw, , brown, fill=brown!40, very thick] (br) at (3.6,0) {};
\node[] at (3.2, 0.4) {Customers};

\node[draw, diamond, fill](v1) at (5.8,0 ) {};
\node[] at (5.8, 0.4) {Vehicles};

\node[] at (7.5, 0.4) {Paths};
\draw[->, dashed, very thick] (7,0) -- (8, 0);

\end{tikzpicture}}
   
   \usetikzlibrary{shapes,snakes,shapes.geometric, patterns, scopes,patterns,intersections,calc}

\begin{tikzpicture}[scale = 0.49]
\useasboundingbox (-6,-5) rectangle (10,10);
\draw (current bounding box.north east) -- (current bounding box.north west) -- (current bounding box.south west) -- (current bounding box.south east) -- cycle;

\node[rectangle, draw] at (-3, 9) {time = 14:45 h};

\node[rectangle, draw, fill, red, very thick] (r) at (1, 1) {};
\node[rectangle, draw, fill, blue, very thick] (b) at (0, -2) {};
\node[rectangle, draw, pattern, brown, very thick] (br) at (3, -1) {};

\node[circle, draw, red, fill=red!40, very thick] (r1) at (1, 6) {\tiny 1};

\node[circle, draw, red, fill=red!40, very thick](r2) at (4, 0) {\tiny 2};

\node[] at (7.4, 0.8) {\small [16:45]$\rightarrow$16:00};
\node[] at (1,7.4) {\small [16:00]$\rightarrow$15:45};
\node[] at (-3.5, 6) {\small [15:10]$\rightarrow$15:15};
\node[] at (2, 2.9) {\small [15:30]$\rightarrow$15:30};
\node[] at (5, -1) {\small [16:45]};

\node[circle, draw, red, fill=red!40, very thick](r3) at (7, 2) {\tiny 3};


\node[circle, draw, blue, fill=blue!40, very thick](b1) at (-3, 5) {\tiny 4};


\node[circle, draw, blue, fill=blue!40, very thick](b2) at (1, 4) {\tiny 5};

\node[draw, diamond, fill](v1) at (-0.5, 1.5) {};
\node[draw, diamond, fill](v2) at (4, -3.5) {};


\draw[->,dashed, very thick] (v1) -- (r);
\draw[->,dashed, very thick] (r) -- (b) ;
\draw[->,dashed, very thick] (b) -- (b1);
\draw[->,dashed, very thick] (b1) -- (b2) ;
\draw[->,dashed, very thick] (b2) -- (r1);
\draw[->,dashed, very thick] (r1) -- (r3);



\end{tikzpicture}
   \usetikzlibrary{shapes,snakes,shapes.geometric, patterns, scopes,patterns,intersections,calc}

\begin{tikzpicture}[scale = 0.49]
\useasboundingbox (-6,-5) rectangle (10,10);
\draw (current bounding box.north east) -- (current bounding box.north west) -- (current bounding box.south west) -- (current bounding box.south east) -- cycle;

\node[rectangle, draw] at (-3, 9) {time = 15:00 h};

\node[rectangle, draw, fill, red, very thick] (r) at (1, 1) {};
\node[rectangle, draw, fill, blue, very thick] (b) at (0, -2) {};
\node[rectangle, draw, fill, brown, very thick] (br) at (3, -1) {};


\node[circle, draw, red, fill=red!40, very thick] (r1) at (1, 6) {\tiny 1};

\node[circle, draw, red, fill=red!40, very thick](r2) at (4, 0) {\tiny 2};

\node[] at (7.4, 0.8) {\small [16:45]$\rightarrow$16:00};
\node[] at (1,7.4) {\small [16:00]$\rightarrow$15:45};
\node[] at (-3.5, 6) {\small [15:10]$\rightarrow$15:15};
\node[] at (2, 2.9) {\small [15:30]$\rightarrow$15:30};


\node[circle, draw, red, fill=red!40, very thick](r3) at (7, 2) {\tiny 3};


\node[circle, draw, blue, fill=blue!40, very thick](b1) at (-3, 5) {\tiny 4};


\node[circle, draw, blue, fill=blue!40, very thick](b2) at (1, 4) {\tiny 5};

\node[] at (-0.3, -0.2) {\small [17:00]$\rightarrow$15:25};

\node[] at (6.5, -0.5) {\small [17:00]$\rightarrow$16:10};

\node[] at (0, -3.5) {\small [17:00]$\rightarrow$15:40};
\node[] at (-3.5, 1) {\small [17:00]$\rightarrow$16:00};


\node[circle, draw, brown, fill=brown!40, very thick](d1) at (0, -1) {\tiny 6};
\node[circle, draw, brown, fill=brown!40, very thick](d2) at (-4, 0) {\tiny 7};
\node[circle, draw, brown, fill=brown!40, very thick](d3) at (-3, -3) {\tiny 8};

\node[draw, diamond, fill](v1) at (-1.4, 3) {};
\node[draw, diamond, fill](v2) at (4, -3.5) {};


\draw[->,dashed, very thick] (v1) -- (b1);
\draw[->,dashed, very thick] (b1) -- (b2) ;
\draw[->,dashed, very thick] (b2) -- (r1);
\draw[->,dashed, very thick] (r1) -- (r3);

\draw[->,dashed, very thick] (v2) -- (br);
\draw[->,dashed, very thick] (br) -- (d1);
\draw[->,dashed, very thick] (d1) -- (d3);
\draw[->,dashed, very thick] (d3) -- (d2);
\draw[->,dashed, very thick] (d2) -- (r);
\draw[->,dashed, very thick] (r) -- (r2);


\end{tikzpicture}
   \caption{Illustration of DPDP-LP solutions associated with Example~1}
   \label{fig:example}
\end{figure}
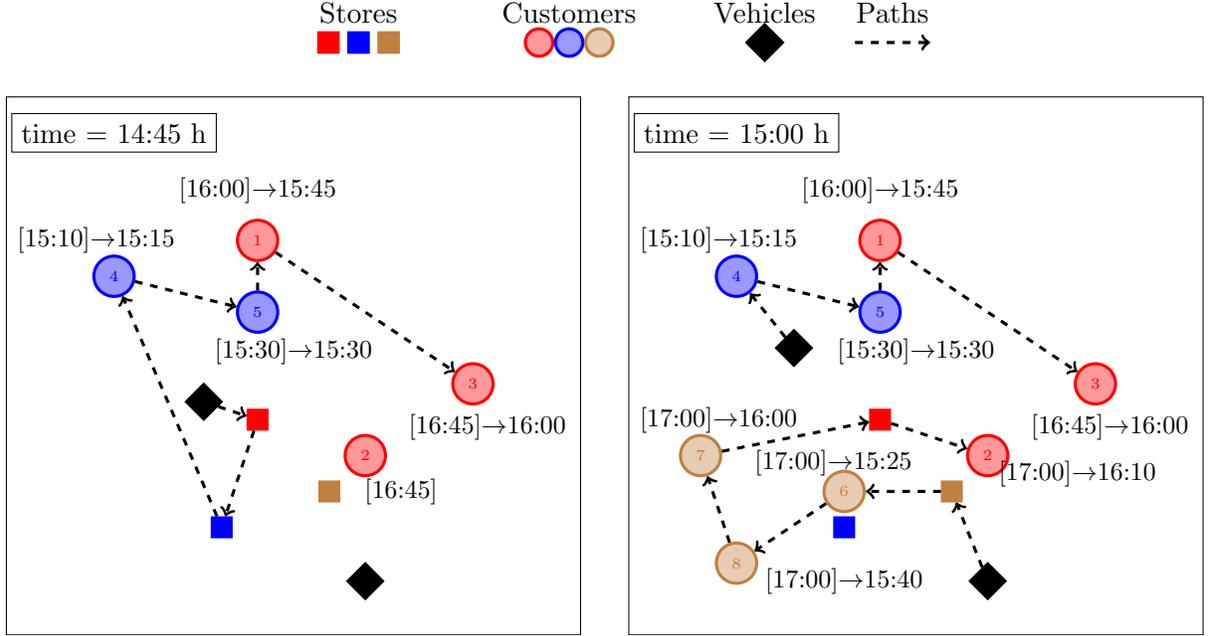
\subsection{Example of the DPDP-LP}

Figure~\ref{fig:example} presents two states of an DPDP-LP instance at times 14:45h and 15:00h. The instance consists of three stores $\cal D = \{1, 2, 3\}$ depicted by the red, blue, and brown square. The two large diamonds indicate vehicles, and customers are represent by nodes in the same color as where their parcel originates from. Each customer has the soft delivery deadline as indicate between brackets, and we depict the actual delivery time of a vehicle with an arrow behind the soft deadline.

In the state at 14:45h, one vehicle has a path assigned to it, starting at the red store, then visiting the blue store, and consecutively visiting the two blue store customers and two red store customers. Observe that the first customer of the tour is visited just after the soft deadline, and a penalty is incurred according to the penalty function $\hat p$. The state at 15:00h shows that the vehicle that was assigned a vehicle path at 14:45h has performed the store visits and is now on its way to the first customer. In the meantime, three new requests appeared in the system each requiring a single product from Store 3 (in brown). A new vehicle path is assigned to vehicle that was idling at 14:45h to visit all the current outstanding requests.

\begin{table}[h!]
    \centering 
    \caption{ Overview of notation used}
    \label{tab:notation}
    \small
    
    \begin{tabular}{l|l}
         \toprule 
         \multicolumn{2}{c}{Problem Description} \\ \midrule
         $\cal T$ & Time horizon, $\cal T = [0, T]$ \\
         $\cal D$ & Set of stores, $\cal D = \{1, \ldots, D\}$ \\
         $\cal V$ & Set of vehicles, $\cal V = \{1, \ldots, V\}$ \\
         $\hat p(t', t^*)$ & penalty of delivering a request at $t^{*}$ with associated delivery deadline $t'$ \\
         \midrule
         \multicolumn{2}{c}{State Variables} \\ \midrule
         $k$ & Time step or decision epoch  \\
         $\cal S$ & State space \\
         $\cal K$ & Set of decision epochs, $\cal K = \{1, \ldots, K\}$ \\
         $S_k$ & State at time step $k$, $S_k \in \cal S$ \\ 
         $t_k$ & Time of decision epoch $k$ \\
         $\cal D_k$ & Vector of store locations of revealed but unassigned requests \\
         $\cal R_k$ & Vector of customer locations of revealed but unassigned requests \\
         $\ell_k$ & Vector of deadlines of revealed but unassigned requests \\
         $a_k$ & Arrival times of vehicles at end of assigned path, $t_k$ if no path assigned \\
         $\hat R_k$ & Customer locations where vehicles becomes idle \\ 
         $n^k$ & Number of unassigned pickup-and-delivery requests at decision epoch $k$ \\
         \midrule
         \multicolumn{2}{c}{Decision Variables} \\ \midrule
         $\cal X(s_k)$ & Space of decision variables \\ 
         $x_k(S_k)$ & Decision at time step $k$ \\
         $C(S_k, x_k)$ & Direct cost of decision $x_k$ \\
         $G^v_k = ((\hat R_{vk} \cup N_k), E_k^v)$ & Graph induced by state $s_k$ for vehicle $v$ \\
         $P^v_k$ & Set of feasible paths of vehicle $v$, $p^v_k \in \cal P^v_k$ \\
         $N^{pvk}$ & Sequence of node visits associated with vehicle path $p^v_k$ \\
         $t^{pvk}$ & Vector of arrival times associated with node visits in vehicle path $p^v_k$ \\
         $\delta_i^{pvk}$ & Binary parameter equalling 1 if node $i$ is visited in path $p^v_k$, 0 otherwise.\\ 
         $y^{pv}_k$ & Binary variable equalling 1 if path $p^v_k$ is assigned to vehicle $v$ at decision epoch $k$ \\
         $c^{pv}_k$ & Direct cost of vehicle path $p^v_k$ at decision epcoh $k$ \\
         $C(S_k, x_k)$ & Direct cost of action $x_k$ in state $S_k$ \\
         \midrule
         \multicolumn{2}{c}{Exogenous Information} \\ \midrule
         $W_{k+1}(S_k, x_k)$ & Exogenous information variable \\
         $\cal D^{\mbox{new}}_{k+1}$ & Vector of store locations associated with new request between $t_k$ and $t_{k+1}$ \\
         $\cal C^{\mbox{new}}_{k+1}$ & Vector of customer locations associated with new requests between $t_k$ and $t_{k+1}$ \\
         $\ell^{\mbox{new}}_{k+1}$ & Vector of deadlines associated with new requests between $t_k$ and $t_{k+1}$ \\
         \midrule \multicolumn{2}{c}{Transition function} \\ \midrule
         $S^M(S_k, x_k, W_{k+1})$ & Transition function \\
         $S^*_k$ & Post-decision state \\
         $\cal R^*_{k}$ & Vector of customer locations in post-decision state \\
         $\cal D^*_{k}$ & Vector of store locations in post-decision state \\ 
         $\ell^*_{k}$ & Vector of deadlines in  post-decision state \\ 
         $a^*_{k}$ & Vector of vehicle arrival times in post-decision state \\ 
         $\hat R*_{k}$ & Vector of vehicle positions in post-decision state \\ 
         \midrule
         \multicolumn{2}{c}{Objective} \\ \midrule 
         $\Pi$ & Set of policies \\
         $X^{\pi}$ & Decision rule, $\pi \in \Pi$ \\ 
         \bottomrule  
    \end{tabular} 
    
\end{table}

\subsection{Markov Decision Process}
\label{sec:mod}

We model the DPDP-LP as a Markov decision process. We first describe the state and decision variables. Then, we define the exogenous information, transition function, and objective function. An overview of all relevant notation is provided in Table~\ref{tab:notation}.

\subsubsection{State Variables.}

At each decision epoch~$k$, we denote the state variable by $S_k, k \in \cal K$. The set of decision epochs $\cal K = \{1, \ldots, K\}$ consists of the time steps at which i) a vehicle finishes an assigned vehicle path and becomes idle again, ii) a new request appears in at the platform, and iii) discrete time steps when vehicles are idle and requests are not assigned. As we consider stochastic dynamic requests, $K$ is a stochastic variable. 


The state variable~$S_k$ is formally described as $S_k = (t_k, \cal D_k, \cal R_k, \ell_k, a_k, \hat R_k)$. We discuss these elements one by one.
\begin{enumerate}
    \item The actual time of the $k^{\mbox{th}}$ decision epoch is denoted by $t_k \in \cal T$.
    \item The information regarding the pickup-and-delivery requests currently revealed but not yet assigned is described  by $n^k$ store locations $\cal D_k \subset \cal D$ and $n^k$ customer locations $\cal R_k$. Thus, $\cal D_k$ and $\cal R_k$ are vectors of size $n^k$. The pickup-and-delivery relation between stores and customers are preserved by their ordering in the vectors $\cal D_k$ and $\cal R_k$, that is, the associated store of customer $R_{ik}$ is $D_{ik}$. With each customer location $R_{ik}$, we denote the associated soft delivery deadline by $\ell_{ik}$. Note, this implies $t_k - (\ell_{ik} - \bar{d})$ equals the time of ordering request $i$ at the platform. 
    \item For each vehicle $v \in \cal V$, we denote with $a_{vk}$ the arrival time at the last customer that is assigned to vehicle $v$. The customer locations at which the vehicles finish their assigned vehicle path are collectively denoted by $\hat{\cal R}_k$. In case a vehicle is idle, $a_{vk} = t_k$.
\end{enumerate}
The initial state $S_0$ is at $k = 0$ with $\ell_k = \cal D_k = \cal R_k = \emptyset$, $a_k = (0, \ldots, 0)$ and $\hat R_k = ((x^0, y^0), \ldots, (x^0, y^0))$.


\subsubsection{Decision Variables.}

At each state~$S_k$ at decision epoch~$k$, we define the space of decision variables by $\cal X(S_k)$. In words, each decision $x_k \in \cal X(S_k)$ describes an assignment of feasible vehicle paths to vehicles currently idling in the system.

For each vehicle $v$, we consider the graph $G^v_k = (\{\hat R_{vk} \cup N_k\}, E_k^v)$ that is induced by the state~$S_k$. The node set~$N_k$ is defined as $N_k = \{D_{1k}, \ldots, D_{n^kk}, R_{1k}, \ldots, R_{n^kk}\}$, and the edge set $E^v_k$ is symmetric and complete. We define the set of feasible vehicle paths $\cal P^v_k$ for vehicle $v$ at epoch~$k$ on the graph $G^v_k$. Each vehicle path $p^v_k \in \cal P^v_k$ describes a series of~$q^{pv}_k$ store and customer locations that we denote by $N^{pv}_k = (n_1^{pvk}, n_2^{pvk}, \ldots, n_{q^{pv}_k}^{pvk})$. The arrival times at each node are denoted by the vector $t^{pv}_k = (t_1^{pvk}, t_2^{pvk}, \ldots, t_{q^{pv}_k}^{pvk})$, here $t_1^{pvk} = t_k$ and the subsequent arrival times can be calculated according to the travel time function $d(\cdot, \cdot)$. Formally, we write $p^v_k = (N^{pv}_k, t^{pv}_k)$. Feasibility of $p^v_k$ implies that
\begin{enumerate}
    \item The first visited node $n_1^{pvk} = \hat R_{vk}$, that is, the vehicle departs from its current position.
    \item The sequence of store and customer locations $N^{pv}_k$ are visited should contain the store location $D_{ik}$ if the associated customer location $R_{ik}$ is visited, and vice versa. That is, we cannot assign a store location without the associated customer location.
    \item If store $D_{ik}$ is included in $N^{pv}_k$, customer $R_{ik}$ is visited after $D_{ik}$.
\end{enumerate}
To encode a decision, we introduce the binary decision variable $y^{pv}_k$ that equals 1 if vehicle path $p^v_k$ is assigned to vehicle~$v$, and that equals 0 otherwise. Furthermore, let~$\delta_i^{pvk}$ be a binary parameter equalling 1 if $i \in N^{pv}_k$, that is, if node~$i$ is part of path $p^v_k$. The decision space $\cal X(S_k)$ is then defined by
\begin{align}
\cal X(S_k) := \left\{ y_k \mid \sum_{ p \in \cal P^v_k} \delta_i^{pvk} y^{pv}_k \leq 1 \ \forall i \in N_k \vee \sum_{p^v_k \in \cal P^v_k} y^{pv}_k \leq 1 \ \forall v \in \cal V \vee  \sum_{p \in \cal P^v_k} p^v_k = 0 \mbox{ if } a_{vk} = t_k \right\}. \label{eq:actionspace}
\end{align}
Thus, the decision space consist of all assignments of feasible vehicle paths to vehicles such that 1) no single store or customer location is assigned more than once, 2) each vehicle is assigned at most one vehicle path, and 3) vehicle paths are only assigned to vehicles that idle at epoch~$k$.

Determining the cost of selecting vehicle path~$p^v_k$ involves evaluating each arrival time in the penalty function, that is, $c^{pv}_k = \sum_{i \in N^{pv}_k}\hat p(\ell_{ik}, t_i^{pvk})$. Note that these costs depend on the shape of $\hat p$, so they are potentially nonconvex and nonlinear. The penalty costs of decision $x_k = y_k$ then equals $C(S_k, x_k) = \sum_{p \in \cal P^v_k}c^{pv}_ky^{vp}_k$.

\subsubsection{Exogenous Information.}

The exogenous information variable $W_{k+1}(S_k, x_k)$ models the arrival of new requests (i.e., in bundles associated with orders) at the system. This induces a set of new store locations to visit $\cal D_{k+1}^{\mbox{new}}$, a set of new customer locations to visit $\cal R_{k+1}^{\mbox{new}}$, and a set of deadlines $\ell^{\mbox{new}}_{k+1}$. Each deadline entry $\ell^{\mbox{new}}_{i,k+1}$ equals $t_{k+1} + \bar{d}$. In case a vehicle becomes idle the aforementioned sets are empty. 

\subsubsection{Transition Function.}
The transition to state~$S_{k+1}$ at time~$t_{k+1}$ depends on the state~$S_{k}$, the action~$x_k$, and the exogenous information $W_{k+1}(S_k, x_k)$ via the transition function $S_{k+1} = S^M(S_k, x_k, W_{k+1})$. First, we consider the (deterministic) transition towards a post-decision state $S_k^*$ that relates to the decision~$x_k$. Second, we consider the (stochastic) transition of the post-decision state $S_k^*$ to $S_{k+1}$, which relates to the exogenous information variable~$W_{k+1}$ modeling potential newly arriving requests.

The transition towards post-decision state $S_k^*$ updates parts of the state variable $S_k$ that are impacted by the assignment of vehicle path~$p^v_k$ according to the assignment variables~$y_k$ in action~$x_k$. The update comprises the following steps:
\begin{enumerate}
    \item We set $a^*_{v,k}$ equal to the arrival time at the last customer location in path $p^v_k$.
    \item For all nodes $n^{pv}_i \in N^{pv}_k$, we remove the associated store $D_{ik}$ from $\cal D_k$ to obtain $\cal D^*_{k}$ if $i \leq n^k$; otherwise, we remove customer $R_{ik}$ from $\cal R_k$ to obtain $\cal R^*_{k}$.
    \item We obtain $\ell^*_{k}$ by removing $\ell_{ik}$ from  $\ell_{k}$.
    \item We set $\hat R^*_{vk}$ equal to the last customer from path $p^v_k$.
\end{enumerate}
Then, the post-decision state $S_{k}^* = (k, \cal R^*_k, \cal D^*_k, \ell^*_k, a^*_k, \hat R^*_k)$. Subsequently, the following steps transition into state $S_{k+1}$:
\begin{enumerate}
    \item We define $\cal R_{k+1} = (\cal R^*_{k}, \cal R^{\textsc{new}}_{k+1})$, $\cal D_{k+1} = (\cal D^*_{k}, \cal D^{\textsc{new}}_{k+1})$, $\ell_{k+1} = (\ell^*_{k}, \ell^{\textsc{new}}_{k+1})$, and $\hat R_{k+1} = \hat R^*_{k}$.
    \item The point of time $t_{k+1}$ is defined as the minimum of the following events:
    \begin{enumerate}
        \item The earliest time a vehicle currently traveling becomes idle, that is, $\min_{v \in \cal V : a_vk \neq t_{k}} a_{vk}$. 
        \item The time at which a new order arrives in the system.
        \item The time $t_{k} + \delta_{\max}$ in case there are unassigned orders and idle vehicles.
    \end{enumerate}
    \item Let $a_{k+1} = a^*_{k}$,  and set $a_{v,k+1} = t_{k+1}$ for idle vehicles, that is, vehicles with no path assigned.
\end{enumerate}
Then, $S_{k+1} = (k+1, \cal R_{k+1}, \cal D_{k+1}, \ell_{k+1}, a_{k+1}, \hat R^*_{k+1})$

\subsubsection{Objective.}

A solution to the DPDP-LP is given by a decision policy $\pi \in \Pi$ and a decision rule $X^\pi: \cal S_k \rightarrow \cal X(S_k)$. Thus, the decision $x_k$ under decision policy $\pi$ is written as $x_k = X^\pi(S_k)$. The objective of the DPDP-LP is then to find a policy that minimizes the expected total penalty cost:
\begin{align}
    \min_{\pi \in \Pi} \mathbb{E}\Big[\sum_{k \in \cal K} C(S_k, X^\pi(S_k))\mid S_0\Big], \label{eq:bell}
\end{align}
where $S_{k+1} = S^M(S_k, X^\pi(S_k), W_{k+1})$. 



\section{Cost Function Approximation (CFA)}
\label{sec:sol}

Solving the Bellman equation~\eqref{eq:bell} to optimality is computationally intractable as we face the curse of dimensionality in the state, action, and transition space. We therefore construct a Cost Function Approximation (CFA) approach to overcome the curse of dimensionality. The CFA uses insights from Approximate Dynamic Programming (ADP) to anticipate future demand and exploits advanced MIP methods to ensure that the exponentially large decision space is explored efficiently. Our CFA is general, in the sense that it does not rely on the combinatorial structure of the DPDP-LP, and thus provides a generic way to combine the power of MIP methods with anticipatory decision making.

In short, our CFA approach adapts the objective and constraints of a set-packing model via a parameterization. This modified set-packing model is then solved using column generation at each decision epoch. Our CFA approach enriches the set-packing model with two parameters~$\alpha$ and~$\beta$ that enables balancing between the actual objective (i.e., the incurred penalty costs), the length of each assigned path (as a proxy for current consolidation opportunities), and the urgency of fulfilling a pickup-and-delivery request. We introduce an \textit{urgency function} that defines for each state variable how urgent an order is. The outcome of the urgency function is scaled with~$\beta$. Next, we scale the length of each assigned path with~$\alpha$. The parameters values for~$\alpha$ and~$\beta$ can be learned offline and are independent of the state. However, the parameters will scale the urgency function and vehicle path lengths, which are both dependent on the state variable. The best-performing value of the parameters can simply be learned by using global optimization tools.

In the following, we refer to the complete solution approach as ``our CFA approach". Below, we first discuss the rationale behind our CFA approach. Then, we introduce the modified set-packing formulation and its parameterization. We then detail how column generation can be used to find high quality solutions and outline how the modification of the set-packing model has a direct impact on the paths being generated via column generation. We end this section by a general outline of our CFA approach.

\subsection{Motivation of our CFA approach}

The CFA approach is motivated by the need to manage the available resources (i.e., vehicles) in such a way that the system not only provides a good service here-and-now, but also anticipates future demand and has resources available to provide good service in the (near) future. In addition, we also need to be able to search the decision space reliably and quickly at each decision epoch. 

We consider a setting where a decision maker can only control the assignment of paths to vehicles. When assigning a path to a vehicle, a decision maker first has to consider the consolidation potential that the vectors of yet unscheduled requests~$(\cal R_k, \cal D_k, \ell_k)$ imposes on the system. The geographical locations of the unscheduled requests will directly determine the ratio between the number of visited locations in a vehicle path and its associated length. If this ratio is low, meaning a long vehicle path with only a few locations to visit, its consolidation potential is said to be low. In case a vehicle paths' ratio is high, meaning relatively many locations are visited in relation to the path's length, its consolidation potential is said to be high. In addition to the consolidation potential, the current time, the soft deadlines of the unscheduled requests, and the penalty function will determine what type of vehicle path is best to assign. A path with a low consolidation ratio could be a good choice if the requests are near their deadline, while a path with a high consolidation ratio could be a bad choice if the deadlines are still far in the future. Yet, a path with a high consolidation ratio could also be a good choice if the path's length is sufficiently short to not impact future resource availability considerably. In other words, we can assign a path to a vehicle here-and-now, but then it should be ``worth it'', meaning, sufficiently many locations should be visited given the length of the assigned path \textit{or} the path fulfills requests of which the deadlines are close. Our CFA approach exactly captures this trade-off. 

Our CFA approach penalizes the length of a vehicle path by scaling it with the parameter~$\alpha$, while it approximates the urgency of a request (i.e., how close it is to its deadline) with a profit scaled via the parameter~$\beta$. These fractions are balanced with the actual penalty function~$\hat p$, so that direct and future costs will be balanced as-best-as-possible. Then, if the urgency profits scaled with~$\beta$ associated with a particular vehicle path outweighs the sum of direct penalty costs and the vehicle path length scaled with~$\alpha$, this vehicle path is said to be approximately profitable and might, depending on set of all possible vehicle paths, be assigned to a vehicle.

\subsection{Modified Set Packing Formulation}

Recall that the we consider vehicle paths $p^v_k \in \cal P^v_k$ with associated penalty costs~$c^{pv}_k$, and that~$\delta^{pvk}_i$ is a parameter equalling 1 if node $i$ is visited in path $p^v_k$. Let $v_{ij}^{pvk}$ be a parameter equalling 1 if location $j \in N^{pv}_k$ is visited directly after location $i \in N^{pv}_k$ in path~$p^v_k$. We then define the modified costs~$\tilde{c}^{pv}_k$ as:

\begin{equation}
    \tilde c^{pv}_k :=  c^{pv}_k + \alpha \sum_{i \in  N^{pv}_k } \sum_{j \in N^{pv}_k} d(i,j)v_{ij}^{pv},
\end{equation}

where $\alpha \geq 0$ balances consolidation opportunities by adding a fraction of the vehicle path's length to the actual objective function that consists of the incurred penalty costs represented by~$c^{pv}_k$.

Besides scaling the length of vehicle paths, our CFA approach also modifies the constraints of the set-packing formulation by considering how close an unscheduled request is to its delivery deadline. For this, we introduce an urgency function $h: \cal T \times \cal R_k \rightarrow \mathbb{R}_{\geq 0}$ mapping the current time $t_k$ and a currently unscheduled customer location $R_{ik}$ to a non-negative real number at decision epoch~$k$.

For each customer location $R_{ik} \in \cal R_k$, we call the associated pickup-and-delivery request ``not urgent'' if~$t_k$ is still far away from its deadline~$\ell_{ik}$, and we call it ``urgent'' if~$t_k$ is close to the request deadline~$\ell_{ik}$. 

Our modified set-packing formulation balances the modified path costs~$\tilde c^{pv}_k$ with the urgency of the requests. To arrive at a minimization problem, urgency profits are incurred if a request is left unassigned at epoch~$k$, and no urgency profits are incurred if the request is assigned. Then, the modified set-packing formulation, which we call the Master Problem at epoch~$k$, asks for solving
\begin{align}
\mathbf{(MP)_k}:= \min \quad & \sum_{v \in \cal V}\sum_{p \in \cal{P}^v_k}\tilde c^{pv}_ky^{pv}_k -  \sum_{r \in R_i} \eta_{rk},  \label{mip:obj}& \\
\mbox{s.t.} \quad & \sum_{p \in \cal{P}^v_k} y^{pv}_k \leq 1 & \forall v \in \cal P^v_k,  \label{mip:vehicleonce}\\
& \sum_{v \in \cal{V}}\sum_{p \in \cal{P}^v_k} \delta^{pv}_ky^{pv}_k \leq 1 & \forall r \in \cal R_k,  \label{mip:locationonce} \\
& \beta \cdot h(t_k, \ell_{rk})(1 - \sum_{v \in \cal{V}}\sum_{p \in \cal{P}^v}\delta^{pvk}_ry^{pv}_k) \leq \eta_{rk} & \forall r \in \cal R_k, \label{mip:beta} \\
& \eta_{rk} \geq 0 & \forall r \in \cal R_k, \label{mip:dom1} \\
&y^{pv}_k \in \{0, 1\} & \forall p^v_k \in \cal P^v_k, v \in \cal V. \label{mip:dom2} 
\end{align}

Objective~\eqref{mip:obj} is the cost-function approximation that balances penalty costs, the scaled routing costs (via~$\alpha$) and the scaled urgency costs (via~$\beta$). Constraints~\eqref{mip:vehicleonce} and Constraints~\eqref{mip:locationonce} ensure that the decision space is modelled correctly, see also the definition of the decision space in epoch~$k$ in Equation~\eqref{eq:actionspace}. Constraints~\eqref{mip:beta} link the modeling variables~$\eta_{rk}$ to the urgency function~$h$. Finally, Constraints~\eqref{mip:dom1}~and~\eqref{mip:dom2} indicate the domain of the decision variables. 

The vehicle path sets~$\cal P^v_k$ are of exponential size and cannot be enumerated. Therefore, we will use the restricted set of paths $\bar{\cal P}^v_k \subseteq \cal P^v_k$ and refer to the~$\mathbf{(MP)_k}$ subject to these restricted sets $\bar{\cal P}^v_k$ as the Restricted Master Problem $\mathbf{(RMP)_k}$. To solve the linear relaxation of~$\mathbf{(RMP)_k}$ to optimality, we need to show that there do not exist vehicle paths in $\cal P^v_k \backslash \bar{\cal P}^v_k$ of negative reduced cost, because if such a path exists it implies that the current solution is suboptimal for $\mathbf{(MP)_k}$. Thus, to solve the linear relaxation of~$\mathbf{(MP)_k}$ to optimality, we iteratively solve the linear relaxation of~$\mathbf{(RMP)_k}$ and check for routes of negative reduced cost that are currently not in the restricted route sets. If we identify routes of negative reduced cost, we include these in the restricted route sets and resolve the linear relaxation of $\mathbf{(RMP)_k}$. We repeat this until we can prove no routes of negative reduced cost exist. This procedure is called column generation \citep[see, e.g.,][]{lubbecke2005selected}

We solve the so-called pricing problem to identify vehicle paths of negative reduced cost. Let $\lambda_{vk}, \pi_{rk}, \mu_{rk} \leq 0$ be the dual variables corresponding to Constraints~\eqref{mip:vehicleonce}-\eqref{mip:beta}, respectively. Thus, the reduced cost $\hat{c}^{pv}_k$ of path $p^v_k \in \cal P^v_k$ is defined as: 
\begin{align}
    \hat{c}^{pv}_k := \tilde c^{pv}_k - \lambda_{vk} - \sum_{r \in \cal R_k} \delta_r^{pvk}[\pi_{rk} - \beta \cdot h(t_k, \ell_{rk})\mu_{rk}].\label{eq:pricing}
\end{align}
The pricing problem then asks for finding the path of minimum reduced cost for each vehicle $v \in \cal V$. Or, equivalently:
\begin{align}
    \bar{p}^v_k = \arg\min_{p^v_k \in \cal P^v_k \backslash \bar{\cal P}^v_k} \hat{c}^{pv}_k
\end{align}
Thus, if $\hat{c}^{\bar p v}_k \geq 0$ for all $v \in \cal V$, the linear relaxation of $\mathbf{(RMP)_k}$ is optimal for the linear relaxation of~$\mathbf{(MP)_k}$.

\subsection{Relation Between $\alpha$, $\beta$ and $\mathbf{(RMP)_k}$}

\begin{figure}[t]
    \centering
    \begin{tikzpicture}
\draw[->, name path=X, very thick] (0,0) -- (6,0) node[anchor=north] {time};
\draw	(0,0) node[anchor=east] {$t_k$}
		(3,0) node[anchor=north] {$\ell_{ik}$}
		(0,0) node[fill=black]{};
\draw   (3,0) node[fill=black]{};

\draw[->] (0,-3) -- (0,6) node[anchor=east] {cost};

\draw[dotted] (3,-3) -- (3,6);

\draw[scale = 1, domain=0:6, smooth, thick, variable = \x, black, dashed, name path=penalty] plot (\x,{0.4*\x});

\draw (6, 2.5) node[anchor =west] {$\hat p$};

\draw[scale = 1, domain=0:6, smooth, thick, red, variable = \x, name path=urgency] plot (\x,{0.5 + 0.5*\x});
\draw[scale = 1, domain=0:6, smooth, thick, red, variable = \x] plot (\x,{0.8 + 0.8*\x});

\draw (6, 3.5) node [red, anchor = west] {$\beta_1 h$};
\draw (6, 5.6) node [red, anchor = west] {$\beta_2 h$};

\draw[scale = 1, domain=0:6, smooth, thick, blue, variable = \x, name path=objective] plot (\x,{-1 * ( 0.5 + 0.1 * \x) });
\draw[scale = 1, domain=0:6, smooth, thick, blue, variable = \x] plot (\x,{-1 * ( 0.8 + 0.4 * \x) });

\draw (6, -3.2) node [blue, anchor = west] {$z (\beta_2)$};
\draw (6, -1.1) node [blue, anchor = west] {$z (\beta_1)$};
\tikzfillbetween[of=penalty and urgency]{red, opacity=0.2} {};
\tikzfillbetween[of=X and objective]{blue, opacity=0.2} {};

\end{tikzpicture}
    \begin{tikzpicture}
\draw[->, name path = X, very thick] (3,0) -- (6,0) node[anchor=north] {time};
\draw[name path = X1, very thick] (0,0) -- (3,0) node[anchor=north] {};

\draw	(0,0) node[anchor=east] {$t_k$}
		(3,0) node[anchor=north] {$\ell_{ik}$}
		(0,0) node[fill=black]{};
\draw   (3,0) node[fill=black]{};

\draw[->] (0,-3) -- (0,6) node[anchor=east] {cost};

\draw[dotted] (3,-3) -- (3,6);

\draw[thick, dashed, name path=penalty1, domain = 0:3] plot(\x, 0);
\draw[thick, dashed] (3,0) -- (3,1);
\draw[scale = 1, domain=3:6, smooth, thick, variable = \x, black, dashed, name path=penalty] plot (\x,{ln(\x - 2) + 1});

\draw (6, 2.5) node[anchor =west] {$\hat p$};

\draw[scale = 1, domain=3:6, smooth, thick, red, variable = \x, name path=urgency] plot (\x,{0.5 + 0.5*\x});

\draw[scale = 1, domain=0:3, smooth, thick, red, variable = \x, name path=urgency1] plot (\x,{0.5 + 0.5*\x});

\draw[scale = 1, domain=0:6, smooth, thick, red, variable = \x] plot (\x,{0.8 + 0.8*\x});

\draw (6, 3.5) node [red, anchor = west] {$\beta_1 h$};
\draw (6, 5.6) node [red, anchor = west] {$\beta_2 h$};

\draw[scale = 1, domain=0:3, smooth, thick, blue, variable = \x] plot (\x,{-1 * ( 0.8 + 0.8 * \x) });
\draw[scale = 1, domain=3:6, smooth, thick, blue, variable = \x] plot (\x,{-1 * ( 0.8 + 0.8 * \x) + ln(\x - 2) + 1 });
\draw[thick, blue] (3,-3.20) -- (3, -2.20);
\draw (6, -3.2) node [blue, anchor = west] {$z (\beta_2)$};

\draw[scale = 1, domain=0:3, smooth, thick, blue, variable = \x, name path = objective1] plot (\x,{-1 * ( 0.5 + 0.5 * \x) });
\draw[scale = 1, domain=3:6, smooth, thick, blue, variable = \x, name path = objective] plot (\x,{-1 * ( 0.5 + 0.5 * \x) + ln(\x - 2) + 1 });
\draw[thick, blue] (3,-2) -- (3, -1);
\draw (6, -1.1) node [blue, anchor = west] {$z (\beta_1)$};

\tikzfillbetween[of=penalty1 and urgency1]{red, opacity=0.2} {};
\tikzfillbetween[of=penalty and urgency]{red, opacity=0.2} {};

\tikzfillbetween[of=X1 and objective1]{blue, opacity=0.2} {};
\tikzfillbetween[of=X and objective]{blue, opacity=0.2} {};

\end{tikzpicture}
    \begin{tikzpicture}
\draw[->, name path = X, very thick] (3,0) -- (6,0) node[anchor=north] {time};
\draw[-, name path = X1, very thick] (0,0) -- (3,0) node[anchor=north] {};

\draw	(0,0) node[anchor=east] {$t_k$}
		(3,0) node[anchor=north] {$\ell_{ik}$}
		(0,0) node[fill=black]{};
\draw   (3,0) node[fill=black]{};

\draw[->] (0,-3) -- (0,6) node[anchor=east] {cost};

\draw[dotted] (3,-3) -- (3,6);

\draw[thick, dashed, name path = penalty1] (0,0) -- (3,0);
\draw[thick, dashed] (3,0) -- (3,1);
\draw[scale = 1, domain=3:6, smooth, thick, variable = \x, black, dashed, name path = penalty] plot (\x,{0.15* (\x - 3)^2 + 1});

\draw (6, 2.5) node[anchor =west] {$\hat p$};

\draw[scale = 1, domain=3:6, smooth, thick, red, variable = \x, name path = urgency] plot (\x,{0.5 + 0.5*\x});
\draw[scale = 1, domain=0:3, smooth, thick, red, variable = \x, name path = urgency1] plot (\x,{0.5 + 0.5*\x});

\draw[scale = 1, domain=0:6, smooth, thick, red, variable = \x] plot (\x,{0.8 + 0.8*\x});

\draw (6, 3.5) node [red, anchor = west] {$\beta_1 h$};
\draw (6, 5.6) node [red, anchor = west] {$\beta_2 h$};

\draw[scale = 1, domain=0:3, smooth, thick, blue, variable = \x] plot (\x,{-1 * ( 0.8 + 0.8 * \x) });
\draw[scale = 1, domain=3:6, smooth, thick, blue, variable = \x] plot (\x,{-1 * ( 0.8 + 0.8 * \x) + 0.15* (\x - 3)^2 + 1 });
\draw[thick, blue] (3,-3.20) -- (3, -2.20);
\draw (6, -3.2) node [blue, anchor = west] {$z (\beta_2)$};

\draw[scale = 1, domain=0:3, smooth, thick, blue, variable = \x, name path = objective1] plot (\x,{-1 * ( 0.5 + 0.5 * \x) });
\draw[scale = 1, domain=3:6, smooth, thick, blue, variable = \x, name path = objective] plot (\x,{-1 * ( 0.5 + 0.5 * \x) + 0.15* (\x - 3)^2 + 1 });
\draw[thick, blue] (3,-2) -- (3, -1);
\draw (6, -1.1) node [blue, anchor = west] {$z (\beta_1)$};

\tikzfillbetween[of=penalty1 and urgency1]{red, opacity=0.2} {};

\tikzfillbetween[of=penalty and urgency]{red, opacity=0.2} {};

\tikzfillbetween[of=X1 and objective1]{blue, opacity=0.2} {};
\tikzfillbetween[of=X and objective]{blue, opacity=0.2} {};

\end{tikzpicture}
    
    \caption{Illustration of the relation between $\beta h(), \hat p()$ and the cost $z(\beta) = \beta h() - \hat p()$ of assigning request  $i$ at decision epoch $t_k$ to some vehicle path, where $\beta_1 < \beta_2$, for three penalty functions $\hat p()$}
    \label{fig:penalty_examples2}
\end{figure}
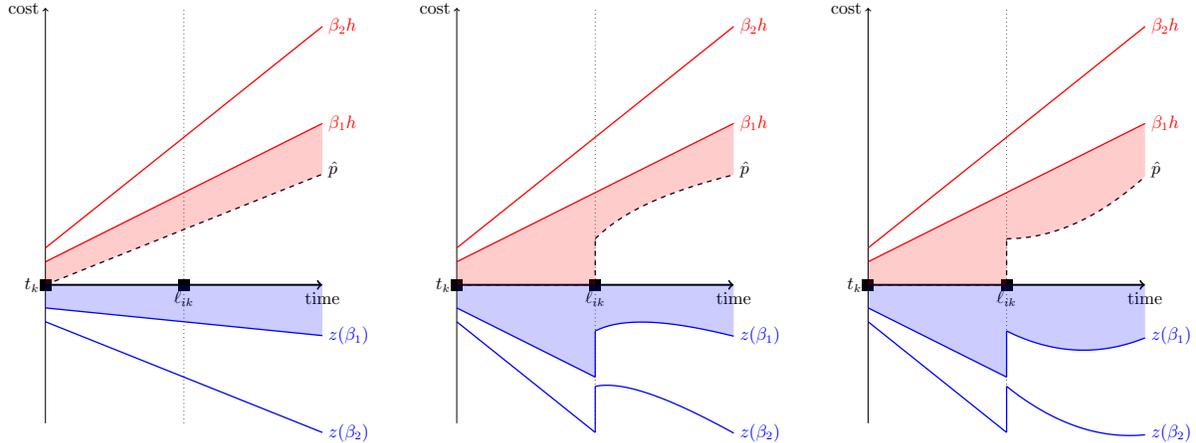

We illustrate the relation between the parameters and the penalty function via the example in Figure~\ref{fig:penalty_examples2}. It shows the same penalty functions as in Figure~\ref{fig:penalty_examples} but now the $x$-axis represents the time $t_k$ of decision epoch $k$. Note that in this figure we do not consider actual delivery times, but focus solely on the time at which a decision is made. The dotted line indicates the minimum penalty incurred at each decision epoch~$t_k$, the red lines~$\beta_1 h$ and~$\beta_2 h$ are two examples of how the urgency function $h$ scales for different values of~$\beta$. In this example, these functions are chosen to be linear with positive intercept. Since the urgency profits are incurred if we do not assign the request at a decision epoch, we plotted the negative of this urgency contribution. Then, the blue lines, indicated by $z(\beta_1)$ and $z(\beta_2)$, plots the \textit{minimum} incurred cost related to request~$i$  if it is contained in an assigned vehicle path at some point along the $x$-axis.

The urgency of a request depends only on the moment its associated vehicle path is assigned to a vehicle, i.e., on~$t_k$ and not on the actual time the vehicle arrives at the customer location. This allows for easy quick evaluation of the reduced cost of vehicle paths and helps to decrease the complexity of finding paths of negative reduced cost via the heuristic approach as explained in the next subsection. The actual penalty costs, however, do depend on the time a vehicle arrives at the customer location. Thus, the blue lines in Figure~\ref{fig:penalty_examples2} indicate the minimum costs of assigning a request (via a vehicle path) to a vehicle at some time along the $x$-axis. Nevertheless, the resulting form of the blue lines clearly show how our parameterization is able to anticipate future demand. The shape of the penalty function affects the shape of the minimum costs incurred for assigned the request via an associated path to a vehicle (i.e., the blue line). In case there is a fixed cost after exceeding the soft deadline (i.e., the middle and right panel in Figure~\ref{fig:penalty_examples2}), this is also reflected in the vehicle path costs. For example, the middle panel shows that our approximation makes assigning a request (via a vehicle path) just before the deadline relatively more attractive then assigning it just after the deadline. Namely, the fixed cost is incurred regardless and for a short moment the request is relatively less urgent and preference should be given to fulfilling requests that still can be delivered timely. 

Finally, we briefly reflect back on the example DPDP-LP solutions illustrated in Figure~\ref{fig:example}. In the state at 14:45h, the first customer is visited after its soft deadline (thus incurring a penalty). If, however, the vehicle would not visit the red store first but directly go to the blue stores and the first two blue customers on the vehicle path, the vehicle could have arrived on time without incurring penalty cost. However, the consequence is that the red store should be visited afterward, resulting in a longer vehicle path. The benefit here is that the red customers (the third and fourth of the assigned vehicle path) are visited earlier, which is good when anticipating future demand. The impact of $\beta$ on the assigned vehicle path is observed in state at 15:00h in Figure~\ref{fig:example}. Clearly, the remaining red customer (that was left unassigned in the state at 14:45h) became urgent enough to be included in the newly assigned vehicle path, so are the new requests.

\subsection{The Dynamic Solution Approach}
We already discussed how column generation can be used generally to solve the linear relaxation of $\mathbf{(MP)_k}$ at each decision epoch $k$. Pricing Problem~\eqref{eq:pricing} is, however, a variant of the Elementary Resource Constrained Shortest Path Problem that is NP-hard in general \citep{dror1994note}. Solving it to optimality using tailored dynamic programming methods \cite[see, e.g.,][]{irnich2005shortest,lozano2016exact, schrotenboer2019branch} is not possible from a computational perspective at each decision epoch $k$. We, therefore, solve Pricing Problem~\eqref{eq:pricing} heuristically to find high quality solutions in short computation times. 

Furthermore, solving the linear relaxation to optimality does not provide an (integer) assignment of paths to vehicles. To find the optimal integer solution to $\mathbf{(MP)_k}$, we embed the (optimal) column generation procedure in branch-and-bound method, yielding a branch-and-price algorithm. However, similar to using a labeling algorithm, the computational burden of a branch-and-price method is too high to be applied at each decision epoch $k$. Instead, we rely on the vehicle paths being generated at the root node only, and do not use column generation in any other nodes of the branch-and-bound-tree. This procedure has already proven its power in other (static) vehicle routing problems \cite[see, e.g.,][]{sonntagtactical,sluijk2021chance}, because the root node optimality gaps are typically small and regularly zero (for smaller instances) for set-packing formulations as the one we use.

In summary, our CFA approach consists of the following steps at each decision epoch~$k$.
\begin{enumerate}
    \item Update the vehicle paths that are currently in $\bar{\cal P}^v_{k - 1}$ based upon the events that happened between $t_{k-1}$ and $t_{k}$.  
    \item Construct the modified set-packing formulation $\mathbf{RMP}_k$ subject to path sets $\bar{\cal P}^v_{k}$, $\alpha$ and $\beta$. 
    \item Repeat the following steps at most $\phi^1$ rounds:
    \begin{itemize}
        \item Solve the linear relaxation of $\mathbf{MP}_k$ , obtain the dual solution $(\bar{\pi}_{rk}, \bar{\lambda}_{vk}, \bar{\mu}_{rk})$. 
        \item For each vehicle $v$ that is currently idle, solve the pricing problem via a stochastic cheapest insertion method. The paths of negative reduced cost that are found via this method are added to the vehicle path sets $\bar{\cal P}^v_k$. If no paths of negative reduced cost are found, go to step~4.
    \end{itemize}
    \item Obtain a new assignment of paths to vehicles by solving $\mathbf{RMP}_k$ to optimality using a MIP solver.
    \item Move to the next decision epoch at time $t_{k+1}$
\end{enumerate}

The stochastic cheapest insertion method is based on the well-known cheapest insertion heuristic and its stochastic variant that reduces complexity. It is typically seen as the insertion method in Adaptive Large Neighborhood Search algorithms \citep[see, e.g.][]{schrotenboer2018coordinating}. In our CFA approach, this heuristic randomly sorts the locations in $\cal D_k$ and $\cal R_k$ and inserts pickup and delivery locations associated with the same request one-by-one on their cheapest position starting with an empty path. Here cheapest position indicates the location of the pickup and delivery in the path that minimizes the associated reduced cost of the path (see Pricing Problem~\eqref{eq:pricing}). This procedure is repeated $\phi^2$ number of times. Each time a path of negative reduced cost is identified, we add it to a buffer that stores the $\phi^3$ paths of most negative reduced cost. Note, in contrast with existing heuristic column generation procedures for vehicle routing problems \citep[see, e.g.,][]{desaulniers2008tabu}, we do not rely on Tabu search methods because preliminary experiments have shown that this method has inferior performance and thereby leads to large variability in the performance of our CFA approach.

\section{Numerical Analysis}
\label{sec:result}

In this section, we discuss the performance of our CFA approach in solving the DPDP-LP compared to various benchmark policies. We start with introducing the base system, including a description of the system parameters that are kept fixed throughout the numerical analysis, in Section~\ref{sec:system}. We then introduce several benchmark policies inspired by prior research and current practice in Section~\ref{sec:benchmark}, including a so-called direct scheduling policy and a policy that limits the length of assigned paths. Section~\ref{sec:perf} presents a comparison between the different policies on key performance indicators. A detailed discussion of the impact of the parameterization on the performance of our CFA approach is given in Section~\ref{sec:perf2}. We end this section by studying the impact of different penalty functions and system characteristics in Section~\ref{sec:sens}.

In the remainder of this section, we consider the following key performance indicators: (i) the expected penalty cost per request, (ii) the expected percentage of pickup-and-delivery requests that are delivered after the soft deadline, (iii) the expected lateness, that is, the expected delivery time minus the deadline under the condition that a pickup-and-delivery request is delivered after its deadline, and (iv) the expected total travel time of the vehicles. For readability, we omit the term `expected' when referring to the key performance indicators.

All the procedures are coded in $C++17$ and use CPLEX~12.8 as MIP solver. The experiments were carried out on an Intel Xeon E5 2680v3 CPU using 4 threads and 8~GB of RAM at most. We rely on enumeration to find the best parameterization for our CFA approach. 

\subsection{Base System}
\label{sec:system}

We consider a base system in which customer locations are uniformly drawn in a $[0,1000] \times [0,1000]$ square. We consider a time horizon of 8 hours during which orders arrive in the system. We introduce an intensity parameter~$\lambda$, denoting the length of the time intervals in which an order is generated with fixed probability 0.2. For example, if $\lambda = 10$ then an  order (both $1-n$ and $n-1$ independently) arrives every 10 minutes in the system with probability 0.2. In the base system, we set $\lambda = 4$ and the size of each order~$n$ is equal to one, meaning that each order includes precisely one pickup-and-delivery request. Each pickup-and-delivery request has a fixed deadline of two hours after it arrived in the system. 

A decision about assigning a vehicle path is made at most every two minutes and at least every five minutes. The former ensures decisions are not made unnecessarily often and the latter ensures requests do not stay unassigned too long without good reason. Throughout the remainder of this section, we denote the time unit in seconds. The travel speed of vehicles equals~0.4 distance units per second. The penalty function equals $\hat{p}(t, r) = \mathbb{I}_{t > \ell^r}(50 + 100(t - \ell^r)/3600)$, that is, there is a fixed penalty of 50 units if the deadline is exceeded and 100 extra units per hour the deadline exceeded. The base system consists of two vehicles $(K = 2)$.

Within our CFA approach, we consider at most 10 pricing rounds at each decision epoch. Each pricing round considers for each vehicle 250 calls to the stochastic cheapest insertion method to solve the pricing problem. In each pricing round, at most 1000 vehicle paths are added; the ones of most negative reduced cost. The final call to the MIP solver solves the RMP to integrality subject to the generated vehicle paths sets and is given at most 20 seconds. However, preliminary experiments show that in most cases root node integrality is achieved, and therefore RMP solves to integrality in negligible time. 

\subsection{Benchmark policies}
\label{sec:benchmark}

To assess the performance of our CFA approach, we consider two benchmark policies that are inspired by current practice at local-to-local shopping platforms and previous research. First, the \textit{Direct Scheduling Policy (DSP)} does not anticipate future demand and assigns requests to idling vehicles immediately at each decision moment. The DSP is obtained from our CFA approach by changing the set-packing constraints of the modified set-packing formulation into set-partitioning constraints and introducing auxiliary variables with high penalty costs to ensure feasibility in case no vehicles are idle. Notice that the DSP still depends on the routing weight parameter~$\alpha$ whereas parameter~$\beta$ is redundant. This policy resembles solution approaches from dynamic pickup-and-delivery problems with restricted solution spaces (e.g., meal-delivery or crowd-sourced delivery). 

Second, the \textit{Limited Length (LIML) policies} assign paths of maximum size~$m$ to idling vehicles at each decision epoch. By setting $\beta = 1,000,000$, we ensure that urgency profits dominate the other costs in our modified set-partitioning formulation. We distinguish among different LIML policies by referring to them as LIML-$m$ policies. In the LIML-1 policy, requests are considered one-by-one in chronological order. In the LIML-$m$ policies with $m > 1$, a policy considers the ``$m$ times the number of idling vehicles'' number of request with earliest deadline, and considers these requests for assignment to vehicles. The LIML policies are inspired by current practice, where relatively short paths are assigned to vehicles in order to easily maintain an overview of the state of the system and free up vehicles to enable a quick response to newly arriving orders. 

\subsection{Comparison between \policy{}, DSP and LIML}
\label{sec:perf}

Here, the performance of our CFA approach is compared to the performance of the benchmark policies DSP and LIML-$m$. 
All results in this section are based on~500 runs of the base system as described in Section~\ref{sec:system}. Although the results are only presented for the base case, an extensive numerical study showed that the insights are representative for a large set of instances. The computation time of a single run takes at most a few seconds, irregardless of the policy being used.

\begin{figure}[t]
    \centering
    \includegraphics[scale=0.70]{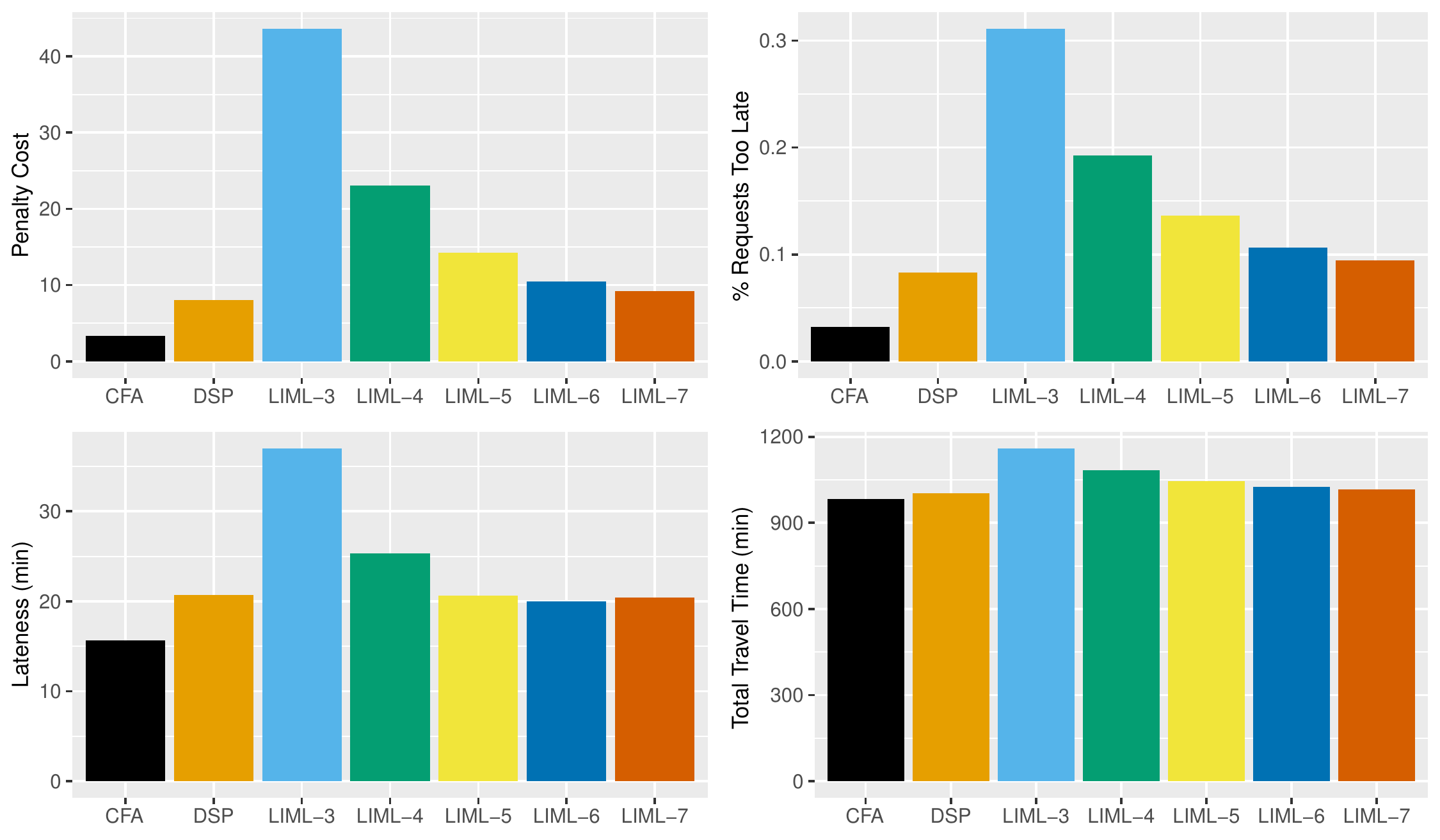}
    \caption{Performance of CFA, DSP, and LIML-$m$ in the base system}
    \label{fig:policyComparison}
\end{figure}

Figure~\ref{fig:policyComparison} presents the penalty cost per request (top left), the \% requests too late delivered (top right), the lateness (bottom left), and the total travel time (bottom right). A few observations stand out. First, our CFA approach clearly outperforms the benchmark policies on all considered performance measures and in particular on the measures that are relevant for customer satisfaction, that is, the \% requests too late and their lateness. Compared to the second-best policy -- the DSP -- our CFA approach reduces the penalty costs per request by 58.3\%, which is mainly due to fewer too late deliveries (8.3\% compared to 3.2\%). Besides that fewer requests are delivered too late, those requests are delivered closer to the deadline too: 16 minutes when using our CFA approach, compared to 21 minutes using the DSP. Second, the total travel times are comparable across the policies except for the LIML-$m$ policies with low values for~$m$, as exemplified by the LIML-3 policy that has a total travel time 17.3\% longer than our CFA approach. Third, although not depicted in Figure~\ref{fig:policyComparison}, the better performance of our CFA approach does not come at the expense of longer delivery times. Indeed, among the orders that are delivered on-time, our CFA approach yields the earliest expected delivery times. Summarizing, our CFA approach outperforms the other policies on customer satisfaction. This is not caused by a increase in travel time, or by an increase of the average delivery times. Thus, using the CFA creates a win-win scenario improving two targets that are usually considered a trade-off, that is, timely delivery and route efficiency.

\begin{figure}[t]
    \centering
    \includegraphics[scale=0.70]{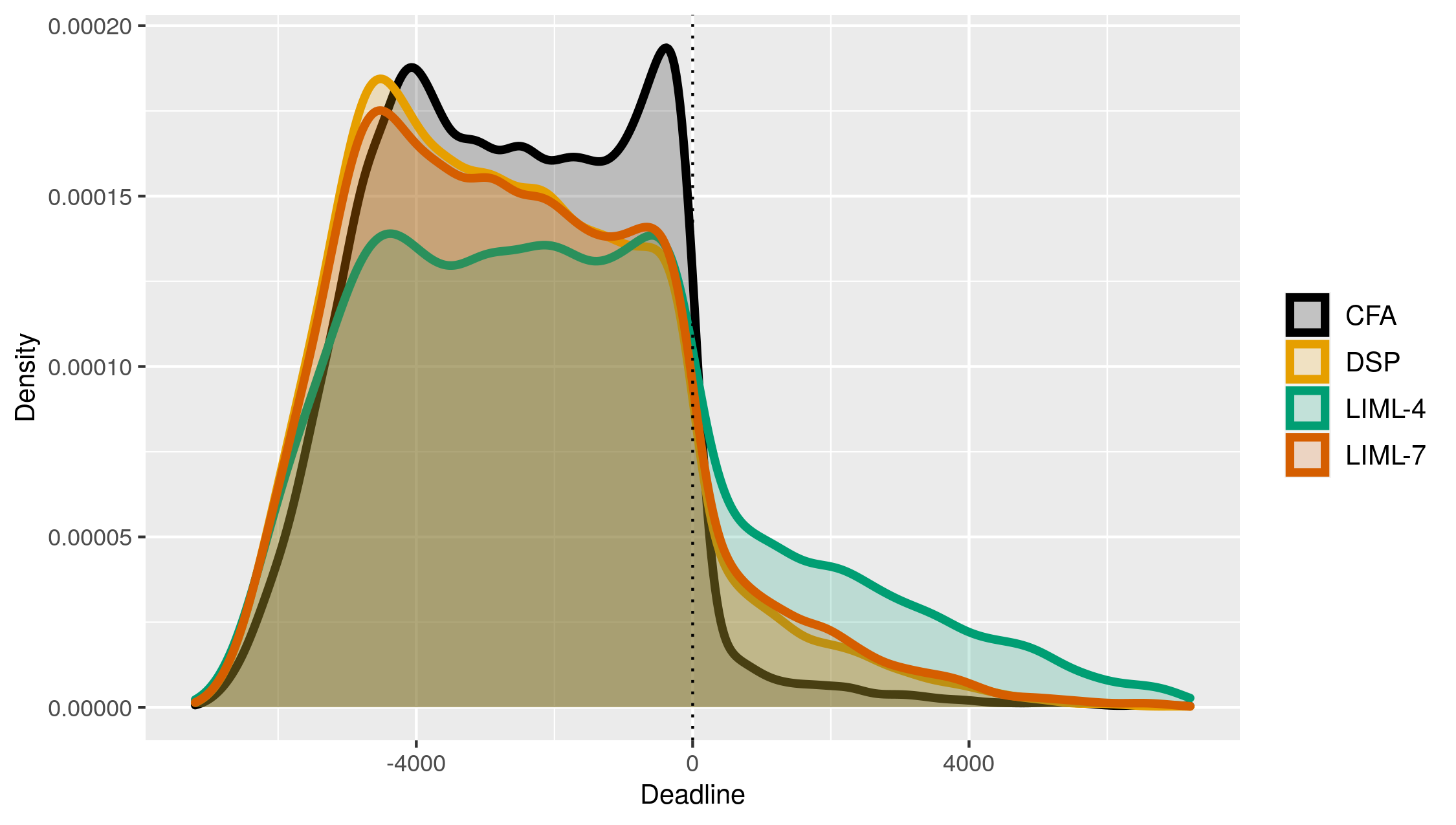}
    \caption{Density of delivery times under CFA, DSP, LIML-4, and LIML-7 in the base system}
    \label{fig:policyComparisonStdev}
\end{figure}

We analyze the delivery times of the requests in more detail in Figure~\ref{fig:policyComparisonStdev}. It shows the empirical density function of the delivery times compared to the deadline associated with each request. We include the performance of our CFA approach, DSP, and LIML-7 -- because these are the best performing policies -- and the performance of LIML-4 to reflect a simple scheduling policy that can be used in practice without large IT investments or advanced routing tools. Positive values indicate deliveries that are delivered too late while negative values indicate on-time deliveries. 

Figure~\ref{fig:policyComparisonStdev} shows that our CFA approach creates more flexibility to realize on-time deliveries in the future by reducing the number of requests that are delivered extremely fast, resulting in the lower probability mass in the upward sloping part -- particularly in relation to DSP and LIML-7. This increased flexibility enables a more efficient use of vehicles overall and allows prioritizing request that are close to their deadline without harming the ability to deliver future requests, leading to a peak in deliveries just before the deadline. Furthermore, the CFA approach has clearly lower density after the deadline. Remarkably, there is also no evidence that the CFA approach sacrifices a few deliveries that are late, as its density is consistently lower then all other policies for positive values. 

\subsection{Performance details of our CFA approach}
\label{sec:perf2}

We continue by studying the impact of~$\alpha$ and~$\beta$ on the performance of our CFA approach. Recall that high values of $\alpha$ emphasize consolidation opportunities by attaching more weight to travel distance of the vehicle paths, whereas high values of $\beta$ emphasize urgency by assigning requests to vehicles relatively long before their deadline.

\begin{figure}[t]
    \centering
    \includegraphics[scale=0.72]{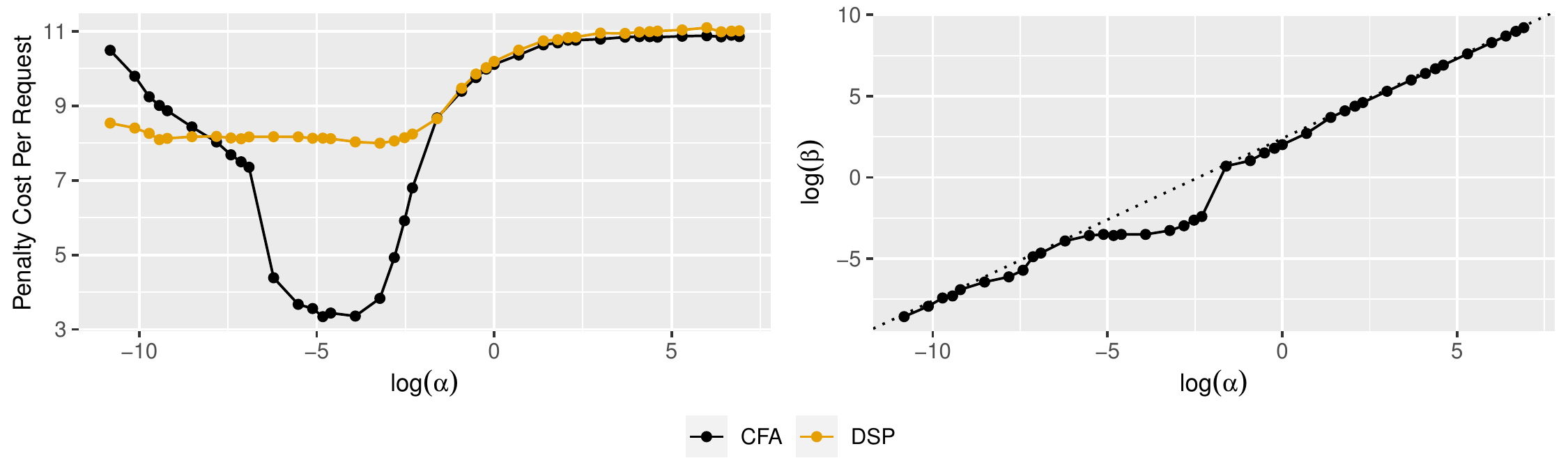}
    \caption{Impact of $\alpha$ on DSP and on CFA with the best associated $\beta$ (left) and relation between $\alpha$ and $\beta$ for CFA (right), on the Base System.}
    \label{fig:alpha1}
\end{figure}

Figure~\ref{fig:alpha1} shows how~$\alpha$ impacts the performance of our CFA approach. The left graph depicts for each value of~$\alpha$ the penalty cost per request of our CFA approach for the best-performing~$\beta$ that is depicted in the right graph. As a benchmark, we also provide the performance of the DSP, but this policy does not depend on~$\beta$ and is therefore omitted in the right graph.

\begin{figure}[b]
    \centering
    \includegraphics[scale=0.72]{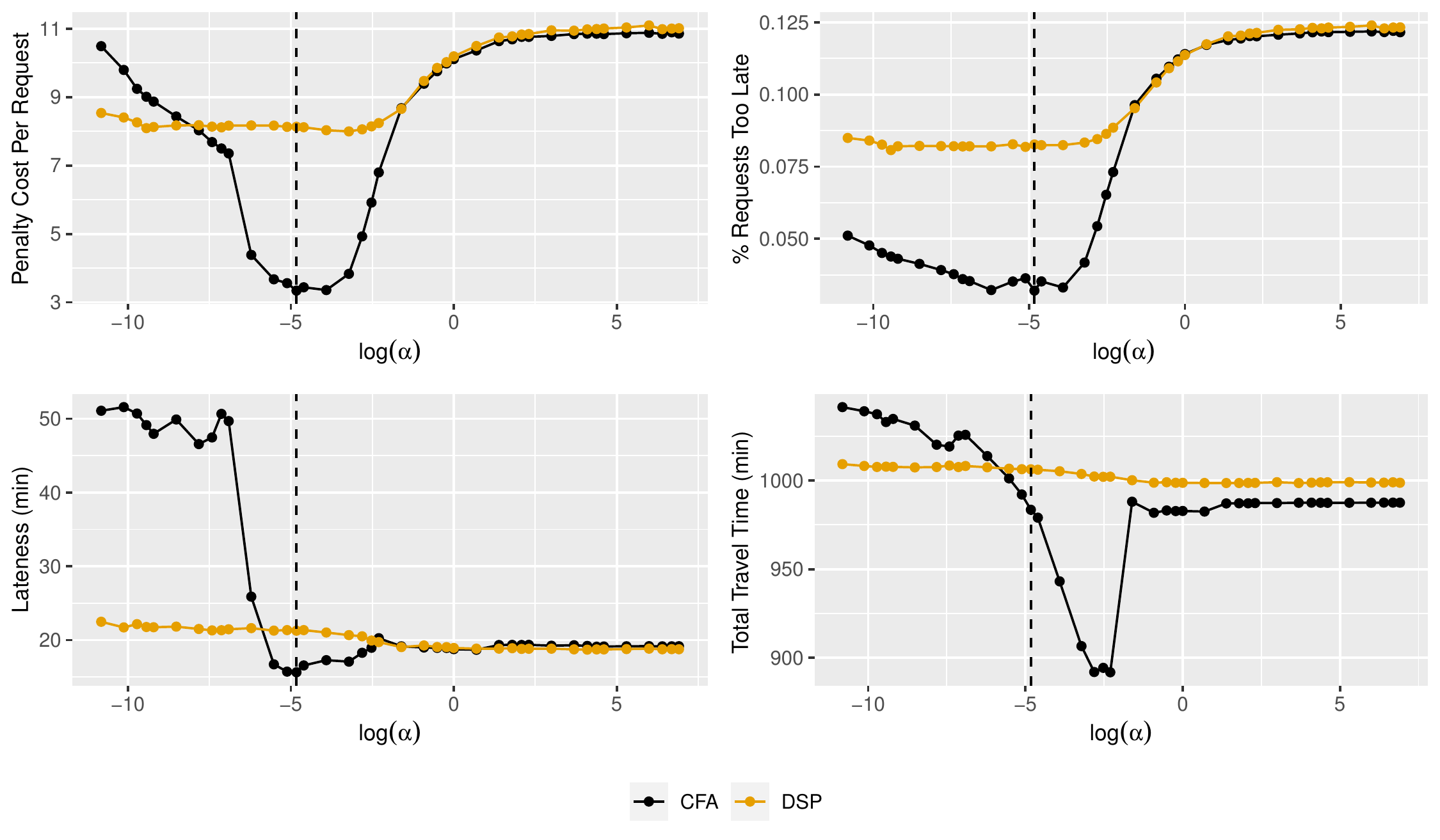}
    \caption{Impact of $\alpha$ on the DSP and CFA with the best associated $\beta$ in the base system}
    \label{fig:alpha2}
\end{figure}
Two observation stand out. First, the penalty cost per request of our CFA approach looks as a smooth function that is first monotonic decreasing and then monotonic increasing for~$\alpha$. If~$\alpha$ is set to a very small number, not enough emphasis is put on the consolidation opportunities compared to the urgency to achieve an efficient policy. The opposite holds for larger values of~$\alpha$, when too much emphasis is put on consolidation opportunities and not on minimizing penalties associated with deadline exceedance. This behaviour is independent of the value of~$\beta$ because each reported point~$\beta$ is chosen so that the penalty cost per request is minimized. Second, the best found values of~$\beta$ depend in a nonlinear way on~$\alpha$. In the right graph, especially around the parameterization that leads to the best performance of our CFA approach, we note that relative smaller values of~$\beta$ are prioritized. Apparently, for the right trade-of between consolidation opportunities (the value of~$\alpha$) and penalty costs, orders should not be assigned too quickly (the value of~$\beta$). It is also interesting to note that for values of~$\alpha$ that are either too small and too large, the best~$\beta$ value depends linearly on~$\alpha$.

\begin{figure}[b!]
    \centering
    \includegraphics[scale=0.75]{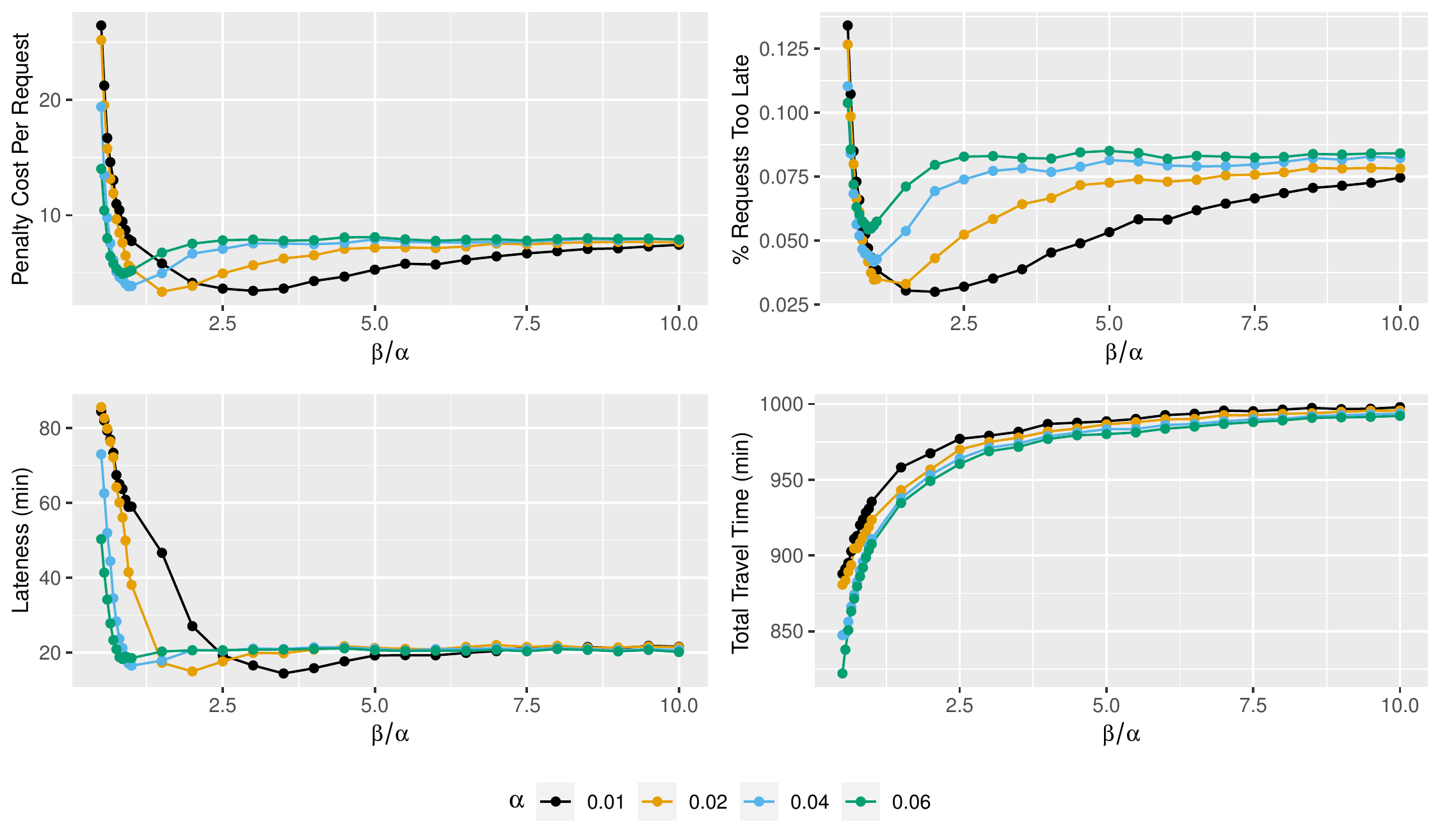}
    \caption{Impact of $\beta$ on CFA performance for various $\alpha$}
    \label{fig:alpha3}
\end{figure}
Figure~\ref{fig:alpha2} presents the relation between~$\alpha$ (with its corresponding cost minimizing~$\beta$) and the performance of our CFA approach and the DSP in more detail. The dashed vertical lines indicate the value of~$\alpha$ for which our CFA approach performs best. A few observations stand out. First, the top left graph shows that for high values of~$\alpha$ the DSP and our CFA approach perform similarly, which can be explained by the relative amount of weight that is put on the vehicle path lengths (as a proxy for consolidation opportunities) compared to the penalty costs and urgency costs (that together are a proxy for the urgency of orders). More importantly, for a small range of~$\alpha$ values close to the optimal~$\alpha$ value, our CFA approach outperforms the DSP by better balancing urgency, penalties, and consolidation opportunities. Notice that the DSP by its design cannot balance these aspects simultaneously. Third, the top right graph shows that the best performing alpha in terms of the penalty cost per request (i.e., the actual objective) also performs well in terms of reducing the \% of too late deliveries, especially in comparison to the DSP. Third, the bottom left graph shows that the lateness is also lower than the DSP. Fourth, the bottom right graph shows that the best-performing policies are not necessarily those that minimize total travel time  distance -- although the total travel time can be reduced considerably, this would come at the expense of higher penalty cost per request.

Figure~\ref{fig:alpha3} further explores the relation between~$\alpha$ and~$\beta$. It shows the same 4 KPI's as in Figure~\ref{fig:alpha2}, but now the $x$-axis depicts the value of~$\beta$ for multiple fixed values of alpha, represented by the coloured lines in the graph. For a given~$\alpha$, the penalty cost per request, the fraction of requests delivered too late, and the lateness all appear to be first monotonic decreasing and then monotonic increasing in~$\beta$. Besides, it becomes apparent that there are unique~$\alpha,~\beta$ combinations that lead to the best-performing CFA policy, as is indicated by the different minima in the top left graph.







\subsection{CFA performance in practical situations}
\label{sec:sens}

We also study the performance of our CFA approach and the various benchmark policies for instances that reflect current practice. We first analyze the performance when the order size increases, that is, when a single order includes multiple products, resulting in multiple pickup-and-delivery requests. We then evaluate how our CFA approach performs for varying penalty functions.

\subsubsection{Varying order sizes.}

To analyze the impact of the order size on the performance of our CFA approach and the benchmark policies, we generated order sizes of the base system randomly between 1 and $n$, where $n$ is varied between 1, 2, and 3. We ensured that instances have the same number of expected pickup-and-delivery requests, thus when $n$ equals 3 there are fewer orders, but of larger size then for $n$ = 1. Table~\ref{tab:resultsN} shows the performance of CFA, DSP, LIML-4, and LIML-7 on the penalty cost per request the \% of requests too late delivered, their lateness, and the total travel time (in minutes) of the vehicles. In addition to the absolute performance values, Table~\ref{tab:resultsN} also presents the relative performance difference (in percentages) compared to our CFA approach.

Several observations stand out. First, for increasing $n$, the LIML-$m$ policies perform relatively worse compared to our CFA approach for all four KPI's considered. Second, the DSP results in slightly higher total travel time compared to our CFA approach, but has more than twice the amount of requests delivered after the deadline, and thus incurs more than twice the expected penalty costs. Third, for increasing $n$, the performance of our CFA approach and the various benchmark policies improve in absolute value, which can be explained by the presence of less orders of smaller sizes implying more requests originating from the same store location at the same time. Finally, we remark that the best-performing parameter values for the policies are similar for varying order sizes. Overall, the performance of our CFA approach is robust for different order sizes. 

\begin{table}[ht!]
    \centering
    \caption{Impact of order size $n$ on the performance of CFA, DSP, LIML-7, and LIML-4}
    \label{tab:resultsN}
    \small
    \begin{tabular}{l|rrr|rrr} \toprule
    & \multicolumn{3}{c}{Absolute values} & \multicolumn{3}{c}{\% difference to CFA} \\
    \cmidrule(lr){2-4} \cmidrule(lr){5-7}
           & $n = 1$       & $n = 2$       & $n = 3$         &  $n = 1$       & $n = 2$       & $n = 3$          \\ \midrule
    \textit{Expected penalty cost}&&&&&& \\
    \ \ CFA    & 3.34 & 1.73& 1.44  &   -        &   -        &          -            \\
    \ \ DSP    & 8.00    & 4.15    & 3.21      & 139.1     & 140.5     & 122.6                \\
    \ \ LIML-7 & 9.23    & 5.48    & 4.75      & 175.9     & 217.6     & 229.2                \\
    \ \ LIML-4 & 23.06   & 14.79   & 13.53     & 589.7     & 756.5     & 837.7                \\[0.25cm]
    \textit{\% Requests delivered after the deadline}&&&&&& \\ 
    \ \ CFA    & 3.20    & 1.80    & 1.39      &       -    &      -     &      -                \\
    \ \ DSP    & 8.33    & 4.39    & 3.38      & 160.2     & 144.2     & 142.9                \\
    \ \ LIML-7 & 9.45    & 5.81    & 4.54      & 195.2     & 223.4     & 226.4                \\
    \ \ LIML-4 & 19.28   & 12.46   & 10.94     & 501.9     & 593.3     & 686.2                \\[0.25cm]
    \textit{Deadline exceedance of too-late delivered requests in minutes} &&&&&&\\ 
    \ \ CFA    & 16      & 9       & 8         &      -     &  -         &       -               \\
    \ \ DSP    & 21      & 15      & 12        & 32.2      & 71.0      & 59.8                 \\
    \ \ LIML-7 & 20      & 14      & 13        & 30.5      & 56.2      & 68.2                 \\
    \ \ LIML-4 & 25      & 18      & 18        & 62.1      & 101.5     & 138.7                \\[0.25cm]
    \textit{Total travel time in minutes} &&&&&&\\ 
    \ \ CFA    & 984     & 888     & 832       & -          &      -     &   -                \\
    \ \ DSP    & 1004    & 912     & 839       & 2.1       & 2.7       & 0.8                  \\
    \ \ LIML-7 & 1017    & 926     & 858       & 3.4       & 4.2       & 3.2                  \\
    \ \ LIML-4 & 1084    & 989     & 918       & 10.2      & 11.3      & 10.4                 \\ \bottomrule
    \end{tabular}
\end{table}
\begin{figure}[h!]
    \centering
    \includegraphics[scale=0.75]{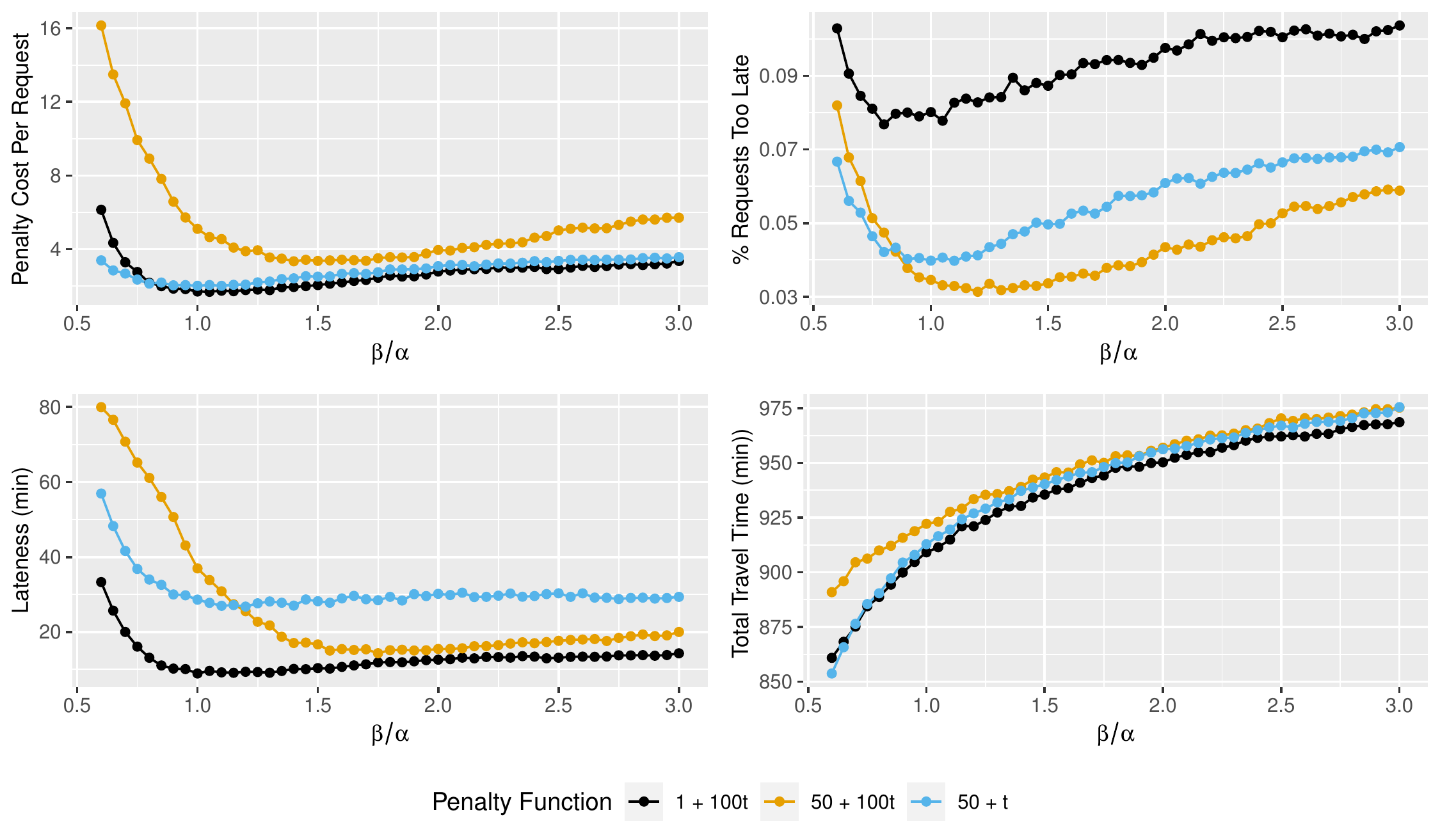}
    \caption{The effect of $\beta$ on the CFA performance for different penalty functions}
    \label{fig:penaltyImpact}
\end{figure}

\subsubsection{Impact of penalty function on parameter values.}

In the previous results, we considered $\hat{p}(t, r) = \mathbb{I}_{t > \ell^r}(50 + 100(t - \ell^r)/3600)$. That is, the penalty function comprised of a fixed cost and a linear variable costs for exceeding a request's soft deadline. Here, we consider alternative penalty functions $\hat{p}_1(t, r) = \mathbb{I}_{t > \ell^r}(50 + 1(t - \ell^r)/3600)$, where the fixed costs is dominating, and $\hat{p}_2(t, r) = \mathbb{I}_{t > \ell^r}(1 + 100(t - \ell^r)/3600)$, where the fixed cost is negligible. In doing so, we study how these different penalty functions impact the parameters of our CFA approach.

Figure~\ref{fig:penaltyImpact} portrays the performance of our CFA approach for the three penalty functions for varying values of $\beta$ and $\alpha = 0.02$, which is representative for all best-performing ranges of $\alpha$. First, we see that for the same value of $\alpha$, the expected penalty cost minimizing values of $\beta$ differ. In other words, the relative importance of urgency differs across the different penalty functions. The relative importance of urgency behaves intuitively. For example, a penalty function with a fixed cost for being late (i.e., the blue line) increases the relative importance of urgency. Figure~\ref{fig:penaltyImpact} also shows that in case no fixed cost is incurred for being late (i.e., the black line), more requests are delivered after their deadline while the average lateness is lower. Interestingly, the lack of a fixed cost for lateness results in a lower total travel time, which can be explained by the fact that a tolerance for (short) delays enables more consolidation opportunities.

\section{Conclusions}
\label{sec:conc}

This paper considers a novel dynamic, stochastic pickup-and-delivery problem in the context of local delivery platforms, called the Dynamic Pickup-and-Delivery Problem for Local Platforms (DPDP-LP). These Local delivery platforms connect local stores with customers by offering instant-delivery of online purchases. Our work is inspired by one of those platforms in the city of Groningen, the Netherlands. We develop a generic solution approach to solve the DPDP-LP. Our approach relies on a two-parameter Cost-Function Approximation (CFA) to anticipate and balance the urgency and consolidation opportunities for future demand and on a column-generation based routing engine to search the extremely large decision-space of pickup-and-delivery problems. The results from this approach are easy to interpret and can be trained completely offline eradicating the need for expensive online scenario evaluations. 

Computational results show the efficiency of our novel approach. Compared to various benchmark policies, our approach improves customer satisfaction, reduces the number of deliveries after their deadlines, reduces their lateness, and it does not increase the total travel time. The approach is also robust under different order sizes, objective functions, and other system characteristics, such as the number of vehicles. A further analysis on the best-performing parameter values of our approach shows that their relation is non-linear, and therefore not trivial to determine without the help of our approach.

The opportunities for further research are numerous. First, our approach is generic and can be applied to other emerging dynamic routing problems too. Especially those that balance consolidation with some form of delivery urgency. Second, local delivery platforms can be further integrated in existing city transportation concepts, such as using parcel lockers, urban consolidation centers, product transfers between vehicles, or the integration of urban freight and public transportation. Third, in practice some orders also originate from, or have a destination, outside the city. Therefore, the integration of local dynamic delivery operations with network design between cities is relevant and an interesting opportunities for further research. Finally, other planning processes related to local delivery platforms, such as workforce planning, are important and could be considered jointly with the dynamic routing, especially in a multiple-day horizon.


\bibliographystyle{informs2014trsc}
\bibliography{ref}


\end{document}